\documentclass[twoside,english,3p]{elsarticle}
\usepackage[T1]{fontenc}
\usepackage{geometry}
\usepackage{amsmath}
\usepackage{amsthm}
\usepackage{stmaryrd}
\usepackage{graphicx}
\usepackage{booktabs}
\usepackage{multirow}
\usepackage{xurl}
\makeatletter
\theoremstyle{plain}
\newtheorem{thm}{\protect\theoremname}
\theoremstyle{boldremark} 
\newtheorem{rem}[thm]{\protect\remarkname}
\ifx\proof\undefined

\providecommand{\proofname}{Proof}
\fi

\renewenvironment{proof}[1][\proofname]{%
	\par\pushQED{\qed}\normalfont%
	\topsep6\p@\@plus6\p@\relax
	\trivlist\item[\hskip\labelsep\bfseries#1\@addpunct{.}]%
	\ignorespaces
}{%
	\popQED\endtrivlist\@endpefalse
}

\journal{Elsevier}

\usepackage{hyperref}
\hypersetup{colorlinks = true, allcolors = blue}

\usepackage[nameinlink]{cleveref}

\crefname{figure}{Fig.}{Figs.}
\crefformat{equation}{Eq.~#2(#1)#3}
\crefformat{section}{Section~#2#1#3}
\AtBeginDocument{%
\let\citet\cite
}

\usepackage[labelfont=bf]{caption}
\@ifundefined{showcaptionsetup}{}{%
 \PassOptionsToPackage{caption=false}{subfig}}
\usepackage{subfig}
\makeatother

\usepackage{babel}
\providecommand{\remarkname}{Remark}
\providecommand{\theoremname}{Theorem}

\begin{document}

\begin{frontmatter}{}

\title{CENN: Conservative energy method based on neural networks with subdomains for solving variational problems involving heterogeneous and complex geometries}

\author[rvt]{Yizheng Wang}

\ead{wang-yz19@mails.tinghua.edu.cn}

\author[rvt]{Jia Sun}

\author[rvt2]{Wei Li}

\author[rvt3]{Zaiyuan Lu}

\author[rvt]{Yinghua Liu\corref{cor1}}

\ead{yhliu@tsinghua.edu.cn}

\cortext[cor1]{Corresponding author}

\address[rvt]{Department of Engineering Mechanics, Tsinghua University, Beijing 100084,
	China}

\address[rvt2]{Department of Mechanical Engineering, Massachusetts Institute of Technology, United States of America}

\address[rvt3]{Faculty of engineering science, KU Leuven, Leuven, 3000, Belgium}

\begin{abstract}
We propose a conservative energy method based on  neural networks with
subdomains for solving variational problems (CENN), where the admissible function satisfying the essential
boundary condition without boundary penalty is constructed by the
radial basis function (RBF), particular solution neural network, and general
neural network. Loss term is the potential energy, optimized based on the principle of minimum potential energy. The loss term at the interfaces has the lower order
derivative compared to the strong form PINN with subdomains. 
The advantage of the  proposed method 
is higher efficiency, more accurate, and less hyperparameters than the strong form PINN with subdomains. 
Another advantage of the proposed method is that it can apply to complex geometries based on the special construction of the admissible function. 
To analyze its performance, the proposed method CENN is used to model representative PDEs,
 the examples include strong discontinuity, singularity, complex boundary, non-linear, 
 and heterogeneous problems. Furthermore, it outperforms other methods when dealing with heterogeneous problems.
\end{abstract}

\begin{graphicalabstract}
	\includegraphics{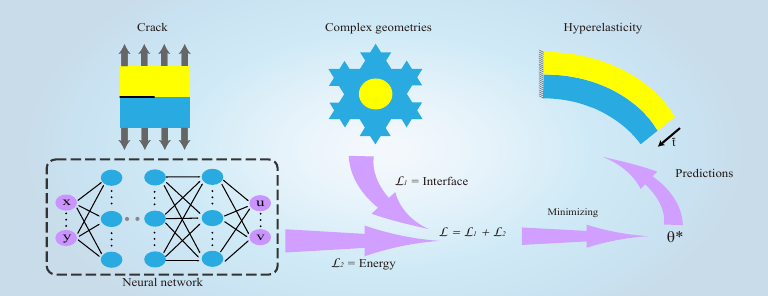}
\end{graphicalabstract}

\begin{keyword}
Physics-informed neural network \sep Deep energy method \sep Domain decomposition
\sep Interface problem \sep Complex geometries \sep Deep neural
network 
\end{keyword}

\end{frontmatter}{}

\section{Introduction}

Many physical phenomena are
modeled by partial differential equations (PDEs). In general, it is difficult
to obtain the analytical solutions of PDEs.
 Hence, various numerical methods are developed to obtain the approximate
solutions in a finite dimensional space. Traditional ways to tackle the solution of PDEs are  the finite element
method (FEM),  the finite difference method,  the finite volume method,
 and mesh-free method \citet{loss_is_minimum_potential_energy}.
  The traditional methods,  especially FEM, are computationally efficient
and accurate for engineering applications.  However, mesh generation of complex boundaries \citet{complex_PINN_a_method_to_construct_admissible_function},  dimensional explosion
and distortion problems for the high-dimensional cases cannot be tackled well by FEM \citet{PINN_review},  reducing
the efficiency exponentially lower. 
In addition, FEM requires the selection of specific basis functions
to construct the approximation function. For example, the singular element
specifically is designed for fracture mechanics singularities \citet{finite_element_book}
in FEM.  Constructing basis functions corresponding
to different elements undoubtedly increases the brainpower cost. 

Artificial intelligence has impacted many fields in the last decade,
 and it is now generally expected that the means to achieve artificial
intelligence can be machine learning, using data-driven optimization.
 Deep learning,  a method in machine learning,  has achieved unprecedented
success in many fields, from  computer vision \citet{ALEXNET},  speech recognition \citet{speech_recognition},
 natural language processing \citet{machine_translation},  and strategy
games \citet{alphago,star_game} to  drug development\citet{alphafold}
.   In the present day, deep learling is ubiquitous, empowering  various
fields. The success of deep learning partially attributes to the
powerful approximation capabilities of the neural network \citet{super_approximation}. It is natural to use the neural network as the approximation
function of PDEs, i.e.,  physics-informed
neural network (PINN) \citet{PINN_original_paper}. 

The idea of using the neural network to solve PDEs can be traced back to the last century \citet{use_neural_network_to_solve_PDE},  but it received little attention due to hardware limitations.  Raissi
et al. propose PINN to deal with the strong form of the PDE,  i.e.,  a weighted
residual method with neural networks as the approximate function \citet{PINN_original_paper}. PINN can apply to many physical systems containing PDEs.
 In \citet{PINN_solid_mechanics}, a framework using PINN was proposed
 for the solution of forward and inverse problems in solid mechanics.
In \citet{PINNfiuld}, PINN is used to  approximate  the  Euler equations  that  model  high-speed  aerodynamic  flows in fluid mechanics.
In \citet{Error_estimates_NS}, the theory of error bound estimation in the incompressible Navier-Stokes equations is proposed for PINN.
 In \citet{PINNbiomechanics}, PINN is used to infer properties of biological materials.
 In \citet{thermochemical_PINN_optimization_way}, PINN is used to
solve the composite material physical system of thermochemical.
  A library named “DeepXDE” using PINN was developed to facilitate the use of scientific machine learning \citet{PINNlibrary}.
 The loss function in PINN is contructed with the PDEs in the domain,  boundary
conditions, and initial conditions by weighted
 sum of the squared errors. Different choices of trial function and test function correspond 
to various numerical methods, eg.  trial function with domain decomposition
\citet{CPINN} and test function with domain decomposition \citet{hp-VPINN}.
 The advantage of PINN strong form is  less dependent on
sampling size in every iteration,  only requiring zero loss at any given coordinate point
\citet{the_comparision_of_strong_and_energy_form}.  Since all PDEs have a weighted residual form,  the strong
form is general and can  solve almost all PDEs.  Another framework of PINN is the deep energy form \citet{the_comparision_of_strong_and_energy_form,admissible_in_PINN_energy_form,PINN_energy_form_to_solve_C0_without_subdomains,PINN_hyperelasticity,deep_ritz,loss_is_minimum_potential_energy,paradeepenergy},
using the physical potential energy as the optimized loss function.
 The advantage of  deep energy method (DEM) is that the physical
 interpretation is stronger than PINN strong form.  In addition,
 DEM requires less hyperparameters than strong form.  By virtue of a smaller derivative order than the strong form, the computational
efficiency and accuracy are higher.  The disadvantage of DEM is dependent on the choice of the integration scheme for the energy
integration,  in order to make the numerical integration of energy as accurate as possible \citet{the_comparision_of_strong_and_energy_form,admissible_in_PINN_energy_form}.
 DEM lacks in generality, not all PDEs have a corresponding energy form \citet{finite_element_book}. 

Most of the current research is focused on the strong form of PINN,
 the research about the energy form is relatively limited. There are too many hyperparameters in the strong
 form of PINN,  especially CPINN (space subdomains) \citet{CPINN} and XPINN (space-time subdomains for arbitrary complex-geometry domains) \citet{XPINN}, so
 we often need to adjust hyperparameters empirically \citet{the_comparision_of_strong_and_energy_form,admissible_in_PINN_energy_form}.
 Although there are currently  valuable research results \citet{NTK_PINN,NTK_to_get_hyperparameter_of_PINN,ill_gradient}
for the selection of hyperparameters, the optimal hyperparameters
still cannot be attained accurately when faced with specific problems. 
The main advantage of domain decomposition of PINN is the flexibility of optimizing all hyperparameters of each neural network separately in each subdomain, so the parallel algorithms can be used to increase the training speed \citet{XPINN_parallel}.
On the other
hand,  the possible displacement field in the DEM is often constructed
by  a penalty factor $\beta$, i.e., $\beta\cdot MSE(u^{pred},\bar{u})$ \citet{deep_ritz},
or  multiplying coordinates, i.e., $x\cdot u^{pred}$ (u=0, when x=0) to satisfy the special geometry essential boundary such as the beam  \citet{loss_is_minimum_potential_energy,PINN_hyperelasticity}. In
\citet{complex_PINN_a_method_to_construct_admissible_function,admissible_in_PINN_energy_form, PINNstrong_form_in_elastodynamics},
the possible displacement field is constructed by distance network
and particular solution network in PINN. Sukumar et al. researched the exact imposition of boundary conditions based on the distance function for the both PINN strong form and energy form \citet{boundary_conditions_distance_functions}. Some of the earlier studies on PINN to construct the admissible function can be traced to the contributions of Lagaris et al. \citet{admissible_earliest_paper}. Although  PINN strong form with subdomains (CPINN)  already exists,   it lacks PINN energy form with
subdomains. 

In this work,  we propose a conservative energy form based on neural
networks with subdomains for solving variational problems (CENN),  where the admissible function is
constructed by the radial basis function (RBF), particular solution
neural network, and general neural network. This method of constructing
the admissible function is suitable for the complexity boundary problem.
To the best of my knowledge, this is  the first attempt to leverage the power of PINN energy
form to the heterogeneous problem, including discontinuity,  singularity,
 high-order tensor, high-order derivative, and  nonlinear PDEs problem. The
advantages of the CENN are multi-fold  : 
\begin{itemize}
\item \textbf{Efficient handling of heterogeneous problems}: Unlike traditional
DEM, CENN can handle the strong discontinuity and derivative discontinuity
problem on the interface by assigning the different neural networks
in each subdomain.  It is worth noting that the interface loss in
CPINN,  has more terms than CENN ,  which will be mentioned in detail
in \Cref{sec:Method}. 
\item \textbf{Hyperparameter fewer}: According to the variational principle, CENN writes the PDEs in the domain as an energy functional and does not consider the hyperparameters of the different PDEs in domains as CPINN to piece PDEs together. In addition, CENN considers fewer interface conditions at the interface than CPINN, which are derived in detail in section 3.2. Therefore, the hyperparameters of CENN are fewer than CPINN. The advantage of fewer hyperparameters is that CENN can reduce the cost of adjusting hyperparameters.
\item \textbf{Complex geometries}: In CENN, the admissible function satisfying the essential boundary condition without boundary penalty is constructed by the RBF, particular solution neural network, and general neural network. CENN can apply to complex geometries based on the special construction of the admissible function.
\item \textbf{Accuracy and efficiency}: Due to the lower derivative in CENN,
 the accuracy and efficiency are higher than the strong form.  In addition,
 the independent part between the subdomains can be implemented
by a parallelization algorithm,  which will further improve efficiency. 
\item \textbf{Less brainpower cost}: CENN benefits from the expressive
power of the neural network.  So we need not construct the approximation
function by designing a basis function,  e.g.  the special elements
in FEM. 
\item \textbf{The ability to solve  distortion problems}:
CENN is a mesh-free method,  so it has the same advantage as the mesh-free
method.  The proposed algorithm is quite effective in distortion. 
\item \textbf{Flexibility of the subdomains configuration}: CENN can divide
the region to different neural networks. The different configurations,
such as the number of hidden layers, activation function,
can be assigned to the different neural networks for the specific problems.
\end{itemize}

The outline of the paper is as follows.  \Cref{sec:Preparatory-knowledge} is a revision of the prerequisite knowledge. It provides a brief introduction
to feed-forward neural network,  PINN, and the DEM.  \Cref{sec:Method} describes the methodology
of the proposed method CENN.  The strategy is explained for constructing the admissible
function based on the RBF and particular network. In \Cref{sec:Result}, some of the representative applications of CENN are presented:
\begin{enumerate}
	\item The crack problem shows that the proposed method can solve
	the strong discontinuity and singularity problem well,  where   data-driven and CPINN are compared with CENN. 
	\item The complex boundary and heterogeneous problem verifies the proposed
	method has the feasibility of complex boundary and derivative discontinuity
	problems,  where various comparisons are made between CPINN, DEM
	and, CENN.
	\item The non-linear hyperelasticity problem with composite materials proves
	the proposed method can deal with non-linear PDEs, where various comparisons
	are made between FEM,  DEM, and CENN.
\end{enumerate}
\Cref{sec:Discussion} shows the discussion. Finally, \Cref{sec:Conclusion} concludes the
study by summarizing the key results of the present work.

\section{Prerequisite knowledge\label{sec:Preparatory-knowledge}}

In this section,  we provide an overview of the feed-forward neural network.
Next,  we introduce the main idea of the PINN. In the end,  we give
an outline of the deep energy method (DEM). 

\subsection{Inroduction to feed-forward neural network}

\begin{figure}[t]
	\begin{centering}
		\subfloat{\centering{}\includegraphics{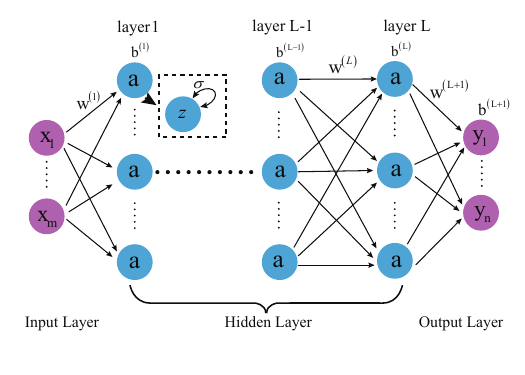}}
		\par\end{centering}
	\caption{Schematic of fully neural network,  the left purple circles  in the neural network
		are the inputs. Blue circles in the neural network are the hidden neurons.
		The right purple circles  in the neural network are the output neurons. $\boldsymbol{x}$ and
		$\boldsymbol{y}$ are the input and output representively. $\boldsymbol{z}$ is the linear output. $\boldsymbol{\sigma}$
		is the activation function acting on $\boldsymbol{z}$.  $\boldsymbol{w}$ and  $\boldsymbol{b}$ are the parameters of the neural network. \label{fig:Schematic-of-NN}}
\end{figure}

The feed-forward neural network is a multiple linear regression with the
activation function aimed to increase the non-linear ability. \Cref{fig:Schematic-of-NN}
shows a schematic diagram of feed-forward neural network. The   feed-forward neural network is given by

\begin{equation}
	\begin{split}\boldsymbol{z}^{(1)} & =\boldsymbol{w}^{(1)}\cdot \boldsymbol{x}+\boldsymbol{b}^{(1)}\\
		\boldsymbol{a}^{(1)} & =\boldsymbol{\sigma}(\boldsymbol{z}^{(1)})\\
		\vdots\\
		\boldsymbol{z}^{(L)} & =\boldsymbol{w}^{(L)}\cdot \boldsymbol{a}^{(L-1)}+\boldsymbol{b}^{(L)}\\
		\boldsymbol{a}^{(L)} & =\boldsymbol{\sigma}(\boldsymbol{z}^{(L)})\\
		\boldsymbol{y} & =\boldsymbol{w}^{(L+1)}\cdot \boldsymbol{a}^{(L)}+\boldsymbol{b}^{(L+1)}
	\end{split}
\end{equation}
where $\boldsymbol{z^{(l)}}$ is the linear transformation of the previous neurons
$\boldsymbol{a^{(l-1)}}$, the layers $1\leq l\leq L$ are the hidden layers.  $\boldsymbol{a^{(l)}}$is
the output of $\boldsymbol{z^{(l)}}$ through the activation function $\sigma$,
and the activation function $\sigma$ is a non-linear function, such
as tanh, sigmoid. In this article,  we uniquely use the
tanh function as the activation function
\begin{equation}
tanh=\frac{e^{x}-e^{-x}}{e^{x}+e^{-x}}.
\end{equation}
where $\boldsymbol{w^{(m)}}$ is the weight between the layer m-1 and m. $\boldsymbol{b^{(m)}}$ is
the bias of the layer m. 

\subsection{Introduction to Physics-Informed Neural Network }

\begin{figure}
	\begin{centering}
		\subfloat{\centering{}\includegraphics{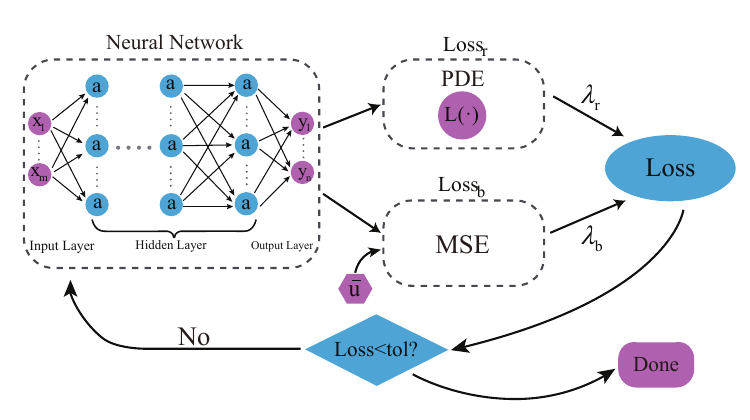}}
		\par\end{centering}
	\caption{Schematic of PINN, the left purple circles  in the neural network are the inputs (temporal and spatial coordinates).
		The right blue circles  in the neural network are the hidden neurons. The right  purple circles 
		in the neural network are the output (interesting field). $\boldsymbol{L(\cdot)}$ is the differential
		operator related to the PDEs. $\bar{\boldsymbol{u}}$ is the given data including
		boundary and initial condition. MSE is the mean square error to
		let the interesting field to satisfy the essential boundary.
		$\lambda_{r}$ and $\lambda_{b}$ are the weight of the residual loss $Loss_r$
		and the boundary loss $Loss_b$ (including initial condition if temporal problem)
		respectively. PINN updates the parameters of neural network until the loss is less than the threshold.
		 \label{fig:Schematic-of-PINN}}
\end{figure}

In this section,  we give the outline for physics-informed neural network (PINN).
The schematic of the PINN is shown in \Cref{fig:Schematic-of-PINN}.
PINN uses  neural networks as the approximation function. PDEs
related to the specific problem are the guide of the loss function
construction. The boundary and initial
conditions are required in advance. The $Loss_{b}$ is the discrepancy
between exact solution and neural
network approximation, which commonly uses MSE as the criterion to satisfy the essential boundary.
Then the $Loss_{r}$ is the difference between the known PDEs
and the neural network approximation, where the differential operators
are regularly constructed by the Automatic Differentiation \citet{automatic_differential}
to get the differential terms in PDEs, which provides exact derivatives while bypassing the computational expense and accuracy issues of symbolic and numercial differentiation \citet{hyper_constitutive_data_driven}. Fortunately, the packages such as Tensorflow
and Pytorch \citet{pytorch_introduction} have the function inherently and conveniently,  so it is convenient to construct
the $Loss_{b}$ related to the PDEs. The whole loss function reads
as :
\begin{equation}
	\mathcal{L}=\frac{\lambda_{r}}{N_{r}}\sum_{i=1}^{N_{r}}|\boldsymbol{L}\boldsymbol{u}(\boldsymbol{x}_{i};\boldsymbol{\theta})|{}^{2}+\frac{\lambda_{b}}{N_{b}}\sum_{i=1}^{N_{b}}||\boldsymbol{u}(\boldsymbol{x}_{i};\theta)-\bar{\boldsymbol{u}}(\boldsymbol{x}_{i})||_{2}^{2}
\end{equation}
where $\boldsymbol{u}(\boldsymbol{x};\boldsymbol{\theta})$ is the prediction of the coordinate point $\boldsymbol{x}$
with the neural network parameters $\boldsymbol{\theta}$. $\boldsymbol{x}$ is the coordinate (spatial
and temporal) which is being used as the input of the neural network;  $\boldsymbol{\theta}$ is
the neural network parameter obtained by the optimization of the loss
function. $\boldsymbol{L}$ is the differential operator,  which can be the linear
or non-linear differential operator. $\boldsymbol{Lu}(\boldsymbol{x};\boldsymbol{\theta})$ is the PDE equation
using the neural network as the interesting field $\boldsymbol{u}$. $\text{\ensuremath{\bar{\boldsymbol{u}}(\boldsymbol{x})}}$is
the boundary or initial condition data known in advance. $N_{r}$
and $N_{b}$ are the number of the residual points and the boundary
data points respectively. $\lambda_{r}$ and $\lambda_{b}$ represent
the corresponding weight of the residual loss $Loss_{r}$ and the
boundary loss $Loss_{b}$ respectively

\begin{equation}
	\begin{split}\mathcal{L}_{r} & =\frac{1}{N_{r}}\sum_{i=1}^{N_{r}}|\boldsymbol{L}\boldsymbol{u}(\boldsymbol{x}_{i};\boldsymbol{\theta})|{}^{2}\\
		\mathcal{L}_{b} & =\frac{1}{N_{b}}\sum_{i=1}^{N_{b}}||\boldsymbol{u}(\boldsymbol{x}_{i};\boldsymbol{\theta})-\bar{\boldsymbol{u}}(\boldsymbol{x}_{i})||_{2}^{2}
	\end{split}
\end{equation}

\subsection{Introduction to deep energy method }

We  introduce the variational principle firstly.
Then we explain how to combine the PINN with the principle of minimum
potential energy. 

If PDEs has variational formulation,  the solution of the strong
form is equal to the variational formulation, i.e., the solution $\boldsymbol{u}$
of strong form makes  $J(\boldsymbol{u})$ stationary with respect to the arbitrary
changes $\delta \boldsymbol{u}$ \citet{finite_element_book}, where $J(\boldsymbol{u})$ represent
the functional and $\boldsymbol{u}$ is a trial function satisfying the essential
boundary condition. If neural network loss is the functional,  we
can get the solution of the strong form by optimization the loss
\begin{equation}
\mathcal{L}=J(\boldsymbol{u})
\end{equation}
\begin{equation}
\boldsymbol{u}^{DEM}(\boldsymbol{x}; \boldsymbol{\theta})=min_{\boldsymbol{\theta}}J(\boldsymbol{u}(\boldsymbol{x};\boldsymbol{\theta}))\approx min_{\boldsymbol{u}}J(\boldsymbol{u}),
\end{equation}
where $\boldsymbol{u}^{DEM}$ is the solution of the DEM. Although the function
space of the neural network is enormous,  there might still be  the approximation
error, i.e., the error between the true function space and the PINN
function space depending on the neural network configuration. $\boldsymbol{\theta}$
is the  parameter of the neural network. It is worth noting that $\boldsymbol{u}$ must
satisfy the essential boundary condition, i.e., Dirichlet boundary
condition. 

\begin{rem}
the trial function $\boldsymbol{u}$ must satisfy the essential boundary condition,
\begin{equation}
	\boldsymbol{u}(\boldsymbol{x})=\bar{\boldsymbol{u}}(\boldsymbol{x}),\boldsymbol{x}\subseteq\varGamma^{\boldsymbol{u}},
\end{equation}
where $\varGamma^{\boldsymbol{u}}$ is the essential boundary,  and $\bar{\boldsymbol{u}}(\boldsymbol{x})$
is the given essential boundary value. In addition,  $\boldsymbol{u}(\boldsymbol{x})\in H^{m}$ if the differential operator has order $2m$, where  $H$ is the Sobolev
space.
\end{rem}

\section{Method\label{sec:Method}}

\subsection{construction of the admissible function\label{subsec:construction-of-the admissible}}

The penalty method can be used to satisfy the essential boundary condition
softly,  but the additional hyperparameter called the penalty factor has
to be considered. The most critical challenge of the penalty method is that the best penalty is not precisely known in advance. Therefore,  it
is necessary to construct the admissible function
\begin{equation}
	\boldsymbol{u}(\boldsymbol{x})=\boldsymbol{u}_{p}(\boldsymbol{x};\boldsymbol{\theta}_{p})+RBF(\boldsymbol{x})\cdot \boldsymbol{u}_{g}(\boldsymbol{x};\boldsymbol{\theta}_{g})\label{eq:admissible function}
\end{equation}
where $\boldsymbol{u}_{p}(\boldsymbol{x};\boldsymbol{\theta}_{p})$ called particular network is a common
shallow network trained on the essential boundary to minimize the
following MSE loss
\begin{equation}
	\mathcal{L}_{p}=\frac{1}{n_{be}}\sum_{i=1}^{n_{be}}(\boldsymbol{u}_{p}(\boldsymbol{x}_{i};\boldsymbol{\theta}_{p})-\bar{\boldsymbol{u}}(\boldsymbol{x}_{i}))^{2}
\end{equation}
$n_{be}$ is the number of the essential boundary points. 

$RBF(\boldsymbol{x})$ called distance network is a radial basis function to give
the nearest distance from $\boldsymbol{x}\in\Omega$ to $\varGamma^{\boldsymbol{u}}$, where
$\Omega$ denotes the domain of interest field. 

Many of the methods approximate the distance function by neural networks \citet{PINNstrong_form_in_elastodynamics, complex_PINN_a_method_to_construct_admissible_function}, which is a universal scheme to approximate the distance function for any complex-shaped boundaries. However, such a treatment of the essential boundary condition is not straightforward since an extra training process is required for neural networks.  Considering that the distance function of the analytical solution is not complex in simple structures, it indeed does not need neural networks with a large amount of iterative calculation to approximate the distance function. The neural networks is suitable for complex boundaries and RBF is suitable for simple boundaries.

The number of RBF
distribution points and the arrangement need to be determined
according to the geometric shape of the essential boundary. Here,  we adopt the commonly used Gaussian function as the basis
function :
\begin{equation}
RBF(\boldsymbol{x})=\sum_{i=1}^{n}w_{i}\phi(|\boldsymbol{x}-\boldsymbol{x}_{i}|),
\end{equation}
where $w_{i}$ is the weight of the center $\boldsymbol{x}_{i}$ in RBF. The value
of $w_{i}$ can be obtained according to the training set \{$\boldsymbol{x}_{i}$,
 $y_{i}$\} , $y_{i}$ is the label of $\boldsymbol{x}_{i}$,  which is the nearest
distance to the essential boundary. $\boldsymbol{x}$ is the coordinate to be evaluated.
The following linear equations  determine $w_{i}$,
\begin{equation}
	\left[\begin{array}{cccc}
		\phi(|\boldsymbol{x}_{1}-\boldsymbol{x}_{1}|) & \phi(|\boldsymbol{x}_{1}-\boldsymbol{x}_{2}|) & \cdots & \phi(|\boldsymbol{x}_{1}-\boldsymbol{x}_{n}|)\\
		\phi(|\boldsymbol{x}_{2}-\boldsymbol{x}_{1}|) & \phi(|\boldsymbol{x}_{2}-\boldsymbol{x}_{2}|) &  & \vdots\\
		\vdots &  & \ddots & \phi(|\boldsymbol{x}_{n-1}-\boldsymbol{x}_{n}|)\\
		\phi(|\boldsymbol{x}_{n}-\boldsymbol{x}_{1}|) & \cdots & \phi(|\boldsymbol{x}_{n}-\boldsymbol{x}_{n-1}|) & \phi(|\boldsymbol{x}_{n}-\boldsymbol{x}_{n}|)
	\end{array}\right]\text{\ensuremath{\left[\begin{array}{c}
				w_{1}\\
				w_{2}\\
				\vdots\\
				w_{n}
			\end{array}\right]}=\ensuremath{\left[\begin{array}{c}
				y_{1}\\
				y_{2}\\
				\vdots\\
				y_{n}
			\end{array}\right]}}\label{eq:RBF=0077E9=009635-1},
\end{equation}
where $\phi(|\boldsymbol{x}-\boldsymbol{x}_{i}|)$ is named after the
radial basis function. Since this function is completely determined
by the distance between $\boldsymbol{x}$ and $\boldsymbol{x}_{i}$ regardless of the direction from
$\boldsymbol{x}$ to $\boldsymbol{x}_{i}$, that is the reason called radial. Since the displacement function needs to be derived, there is a higher requirement for the continuity of the distance function. However, the analytical distance function is C0 continuous (at the medial axis of the domina, the displacement function is shape) \citet{boundary_conditions_distance_functions}, so the analytical displacement function is not suitable to use in PINN theoretically. Therefore, we choose Gaussian kernel function as the basis function, which ensures the continuity of higher order. Note that the displacement function does not need to be exactly the same as the analytical solution, but only if it is equal to zero on the essential boundary, which guarantees that the inhomogeneous boundary condition is satisfied. Gaussian kernel is:
\begin{equation}
\phi(r)=exp(-\gamma r^{2})
\end{equation}
In this article, we adopt $\gamma=0.5$ as the shape parameter empirically, because the optimal shape parameter is often located at the small scope. At present, the basis function and shape parameter in RBF mainly depend on experience, and there is no good theoretical guidance \citet{Choosing_basis_functions_and_shape_parameters, comparision_RBF}. It is appropriate that the placement of the RBF points reflect the given value surface as well as possible. A good placement of the center points improves the accuracy of the approximation function. We use a uniform distribution way due to the simplicity (other collocation ways: Halton and epsilon collocation) \citet{comparision_RBF}. The method of box counting can be adopted to determine the number of collocation points. Box counting is a kind of a method to determine the complexity of space \citet{box_counting}. For example, we can use the box counting to determine the complexity of the essential boundary conditions. The number of collocation points can be proportional to the box counting, because the large box counting means a more complex the boundary (the more points are needed). The number
of center points  $\{\boldsymbol{x}_{i},y_{i}\}$ commonly often does not need to take much to construct
the distance network, as shown in the later numerical experiment.
It is not difficult to prove that the radial basis matrix on the LHS
of \Cref{eq:RBF=0077E9=009635-1} is invertible if the coordinates of the center points
 are different. RBF method can accurately satisfy the label
value $y_{i}$ of the fixed point $\boldsymbol{x}_{i}$. The label value $y_{i}$
can be obtained in advance using kdtree \citet{complex_PINN_a_method_to_construct_admissible_function}
or other methods. The accuracy can be improved by increasing the
RBF center points, but it will increase the amount of calculation. 

\begin{figure}
	\begin{centering}
		\subfloat{\centering{}\includegraphics[scale=0.625]{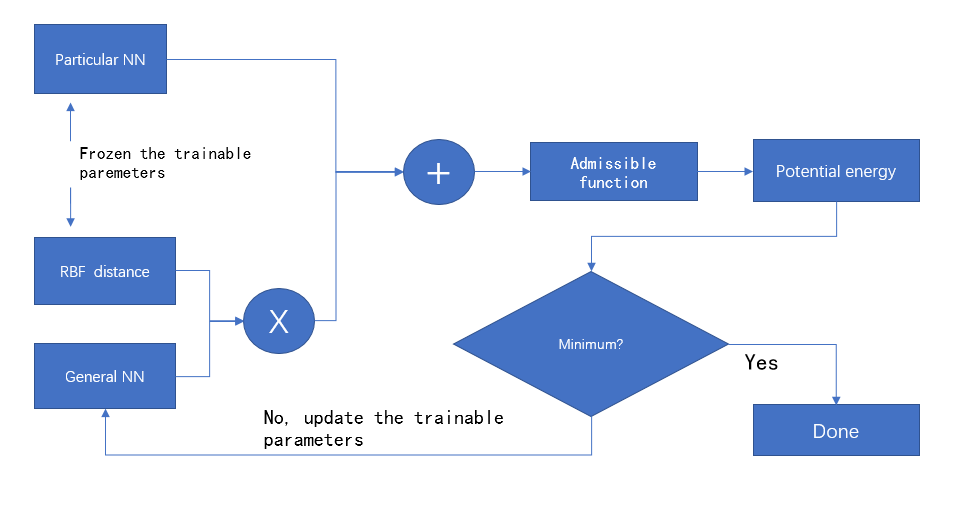}}
		\par\end{centering}
	\caption{The overall process of PINN energy form.\label{fig:Schematic-of-DEM}}
\end{figure}

$\boldsymbol{u}_{g}(\boldsymbol{x};\boldsymbol{\theta}_{g})$   in \Cref{eq:admissible function}  called general network is the regulator to
let the functional $J(\boldsymbol{u})$ minimum. $\boldsymbol{u}_{g}(\boldsymbol{x};\boldsymbol{\theta}_{g})$ is also
a neural network,  whose architecture and activation can be determined
according to the specific problem. Note that the form of \Cref{eq:admissible function}
ensures that $\boldsymbol{u}(\boldsymbol{x})$ satisfies the essential boundary condition if the
particular and distance network have been trained successfully. We
can see the later numerical experiments show it is easy to train the
particular and the distance network successfully. It is important
that the admissible function can solve the complex boundary problem,
 which will be discussed in detail in \Cref{sec:non homogeneous boundary problem}.
 PINN energy method must meet the admissible function before using the principle of minimum potential energy to optimize the potential energy (loss function), as shown in \Cref{fig:Schematic-of-DEM}.

\begin{rem}
The parameter of the particular and the distance network,
which should be trained before the general network in advance,  must
be fixed (untrainable
) when training the general network. The order of particular and distance network does not matter.
\end{rem}

\subsection{CENN}

\begin{figure}
\begin{centering}
\subfloat{\centering{}\includegraphics{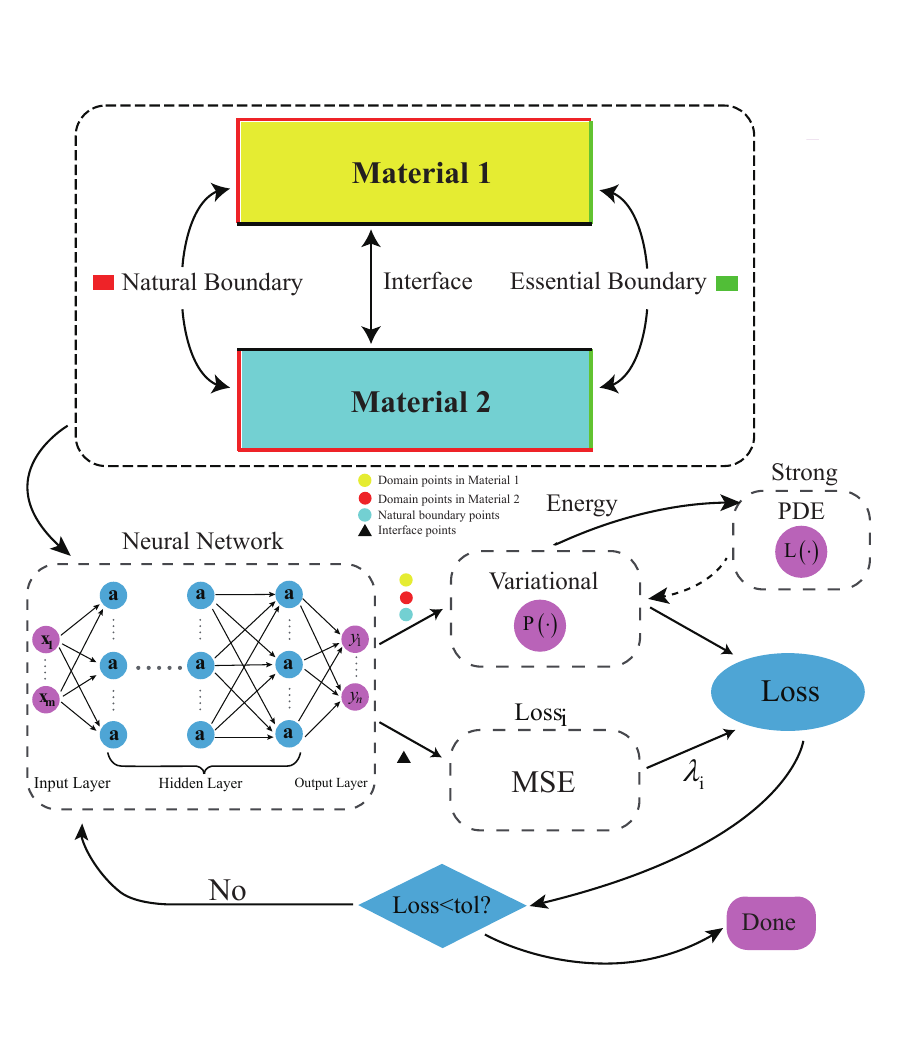}}
\par\end{centering}
\caption{Schematic of CENN : The green lines in the upper picture are the essential
boundary points, which is used to train the admissible function.
The red lines in the upper picture are the natural boundary points. The yellow
and blue region in the upper picture are the domain points corresponding
to the different material. The interface points are the black line
in the upper picture. $\boldsymbol{L}(\cdot)$ is the differential operator related
to the PDE strong form. $\boldsymbol{P}(\cdot)$ is the differential operator related
to the energy  form, whose order of derivative is lower than the strong
form. $\lambda_{i}$ is the weight of the interface term. The input is the coordinates. The output is the field of interest. The loss function is the energy functional. The energy functionals of different materials are connected through the interface loss.\label{fig:Schematic-of-CENN}}
\end{figure}

In this section,  we introduce our proposed method CENN (conservative energy neural network for solving variational problems), and the difference with CPINN (Conservative physics-informed neural networks). Notes, CENN is only used for the PDEs that have variational formulation, but CPINN can be used for almost any PDEs. Note that CENN is based on the principle of the minimum potential energy, so it only can solve the static problem.

\Cref{fig:Schematic-of-CENN} shows the schematic of CENN. CENN is
a deep energy method with subdomains and an admissible function.
The way to construct the admissible function is illustrated in \Cref{subsec:construction-of-the admissible}.
The solution is obtained by optimizing the energy function. The energy
function adds the additional term called interface loss,  which ensures
the conservative physic law.  

We consider the heterogeneous PDEs governing equation :
\begin{equation}
	\begin{split}\boldsymbol{L}(\boldsymbol{u};a_{i})+\boldsymbol{f}_{i}=0, & \boldsymbol{x}\subseteq\Omega_{i}\\
		\boldsymbol{u}(\boldsymbol{x})=\boldsymbol{\bar{u}}(\boldsymbol{x}), & \boldsymbol{x}\subseteq\varGamma^{\boldsymbol{u}}\\
		\boldsymbol{B}(\boldsymbol{u};a_{i})=\boldsymbol{\bar{t}}, & \boldsymbol{x}\subseteq\varGamma^{\boldsymbol{t}_{i}}\\
		\boldsymbol{B}(\boldsymbol{u};a^{+})=\boldsymbol{B}(\boldsymbol{u};a^{-}), & \boldsymbol{x}\subseteq\varGamma^{inter}\\
		\boldsymbol{u}^{+}(x)=\boldsymbol{u}^{-}(x), & \boldsymbol{x}\subseteq\varGamma^{inter},
	\end{split}
\end{equation}
where $\boldsymbol{L}$ is the differential operator. $\boldsymbol{u}$
is an interesting tensor field. $a_{i}$ is the coefficient of $\boldsymbol{L}(\boldsymbol{u};a_{i})$
in the different region $\Omega_{i}$. $\boldsymbol{u}(\boldsymbol{x})=\boldsymbol{\bar{u}}(\boldsymbol{x})$
is the essential boundary condition, and $\boldsymbol{\bar{u}}$ is
the given value of the essential boundary condition (Dirichlet boundary
condition). $\boldsymbol{B}(\boldsymbol{u};a_{i})=\bar{\boldsymbol{t}}$
is the natural boundary condition, and \textbf{$\bar{\boldsymbol{t}}$}
is the given value of the natural boundary condition (Neumann boundary
condition). $\Omega$, $\varGamma^{\boldsymbol{u}}$,  $\varGamma^{\boldsymbol{t}}$ and $\varGamma^{inter}$ are
the domain,  essential boundary,  natural boundary and interface respectively.
The subscript $i$ represent the different region. The $+$ and $-$
represent the regions on either side of the interface. We can get
the corresponding integral Galerkin form,

\begin{equation}
	\begin{split}\int_{\Omega_{i}}(\boldsymbol{L}(\boldsymbol{u};a_{i})+\boldsymbol{f}_{i})\cdot\delta\boldsymbol{u}d\Omega+\int_{\varGamma^{t_{i}}}(\boldsymbol{B}(\boldsymbol{u};a_{i})-\bar{t})\cdot\delta\boldsymbol{u}d\varGamma\\
		\int_{\varGamma^{inter}}(\boldsymbol{B}(\boldsymbol{u};a^{+})-\boldsymbol{B}(\boldsymbol{u};a^{-}))\cdot\delta\boldsymbol{u}d\varGamma+\int_{\varGamma^{inter}}(\boldsymbol{u}^{+}-\boldsymbol{u}^{-})\cdot\delta\boldsymbol{u}d\varGamma & =0,
	\end{split}
\end{equation}

If the operator $\boldsymbol{L}$ satisfies
\begin{align}
	\int_{\Omega}\boldsymbol{L}(\boldsymbol{u})\cdot\boldsymbol{v}d\Omega & =\int_{\Omega}\boldsymbol{u}\cdot\boldsymbol{L}^{*}(\boldsymbol{v})d\Omega+b.t.(\boldsymbol{u},\boldsymbol{v})\\
	\boldsymbol{L}(a\boldsymbol{u_{1}}+b\boldsymbol{u_{2}}) & =aL(\boldsymbol{u_{1}})+bL(\boldsymbol{u_{2}}),
\end{align}
we call the differential operator is self-adjointness and linear
operator \citet{finite_element_book}. We assume $\boldsymbol{L}$ is the self-adjointness
and linear operator. We use integration by parts and obtain
\begin{equation}
	\delta J=0
\end{equation}
where 

\begin{equation}
	J=\int_{\Omega_{i}}[\frac{1}{2}\boldsymbol{u}\cdot\boldsymbol{L}(\boldsymbol{u};a_{i})+\boldsymbol{u}\cdot\boldsymbol{f}]d\Omega+\int_{\varGamma^{inter}}(\boldsymbol{u}^{+}-\boldsymbol{u}^{-})^{2}d\varGamma+b.t.(\boldsymbol{u})
\end{equation}

\begin{rem}
	If the differential operator is linear and self-adjoint, and the highest
	differential operator order is even, then strong form must have the
	corresponding minimum potential energy formualtion, i.e., the energy
	formulation is extreme value problem \citet{finite_element_book}.
	In other words, all strong form have a weak form, but only the linear,
	self-adjoint and even order differiential operator problem have the
	corresponding minumum potential energy formulation.
\end{rem}

For the sake of simplicity,  we consider
the specific PDEs,

\begin{equation}
	\begin{split}-a_{i}\triangle(u(\boldsymbol{x}))=0, & \boldsymbol{x}\subseteq\Omega_{i}\\
		u(\boldsymbol{x})=\bar{u}(\boldsymbol{x}), & x\subseteq\varGamma^{u}\\
		\boldsymbol{n}\cdot a_{i}(\nabla u(\boldsymbol{x}))=\bar{t}(\boldsymbol{x}), & \boldsymbol{x}\subseteq\varGamma^{t_{i}}\\
		\boldsymbol{n}\cdot(a^{+}(\nabla u(\boldsymbol{x}))=\boldsymbol{n}\cdot(a^{-}(\nabla u(\boldsymbol{x})), & \boldsymbol{x}\subseteq\varGamma^{inter}\\
		u^{+}(\boldsymbol{x})=u^{-}(\boldsymbol{x}), & \boldsymbol{x}\subseteq\varGamma^{inter}
	\end{split}
\end{equation}
where $\triangle$ is the Laplace operator.  The interesting field
$u$ is a scalar or tensor,  for the sake of simplicity,  we consider
$u$ as a scalar. $\bar{u}$ is the given value of the essential boundary
condition (Dirichlet boundary condition). $\bar{t}$ is the natural
boundary condition (Neumann boundary condition). $\Omega$, $\varGamma^{u}$,
 $\varGamma^{t}$ and $\varGamma^{inter}$ is the domain,  essential
bound,  natural boundary and interface respectively. $a$ is the coefficients
of PDE and constant in the different regions.  The subscript $i$
represent the different region. $n$ is the normal direction of the
Neumann boundary and the interface. For the sake of simplicity,  we
analyze two different regions denoted as $+$ and $-$. The variation
form is constructed by an integral Galerkin form

\begin{equation}
	\begin{split}\int_{\Omega^{+}}-a^{+}(\triangle u)\delta ud\Omega+\int_{\Omega^{-}}-a^{-}(\triangle u)\delta ud\Omega\\
		+\int_{\varGamma^{t+}}(\boldsymbol{n}\cdot a^{+}(\nabla u)-\bar{t})\delta ud\varGamma+\int_{\varGamma^{t-}}(\boldsymbol{n}\cdot a^{-}(\nabla u)-\bar{t})\delta ud\varGamma\\
		+\int_{\varGamma^{inter}}\boldsymbol{n}\cdot(a^{+}(\nabla u)-a^{-}(\nabla u))\delta ud\varGamma+\int_{\varGamma^{inter}}(u^{+}-u^{-})\delta ud\varGamma & =0
	\end{split}
	\label{eq:variation form}
\end{equation}
where $\boldsymbol{u}$ statisfies the essential boundary in advance. We integrate
the  \Cref{eq:variation form} by parts,  and we can get the weak form
with the interface term
\begin{equation}
	\begin{split}\int_{\Omega^{+}}a^{+}(\nabla u)\cdot(\nabla\delta u)d\Omega+\int_{\Omega^{-}}a^{-}(\nabla u)\cdot(\nabla\delta u)d\Omega\\
		-\int_{\varGamma^{t+}+\Pi^{+}}a^{+}\boldsymbol{n}\cdot(\nabla u)\delta ud\varGamma-\int_{\varGamma^{t-}+\Pi^{-}}a^{-}\boldsymbol{n}\cdot(\nabla u)\delta ud\varGamma\\
		+\int_{\varGamma^{t+}}(\boldsymbol{n}\cdot a^{+}(\nabla u)-\bar{t})\delta ud\varGamma+\int_{\varGamma^{t-}}(\boldsymbol{n}\cdot a^{-}(\nabla u)-\bar{t})\delta ud\varGamma\\
		+\int_{\varGamma^{inter}}\boldsymbol{n}\cdot(a^{+}(\nabla u)-a^{-}(\nabla u))\delta ud\varGamma+\int_{\varGamma^{inter}}(u^{+}-u^{-})\delta ud\varGamma & =0
	\end{split}
\end{equation}
where $\Omega^{+}$ and $\Omega^{-}$ are the different regions of the
interesting field. $\varGamma^{+}$ and $\varGamma^{-}$ are Newmann
boundaries of the different region without the interface. $\Pi^{+}$
and $\Pi^{-}$are the interface $\varGamma^{inter}$ of the different
region,  overlapping each other, and $\boldsymbol{n}$ is the normal direction
of the region $\Omega^{+}$ on the interface, i.e., $\boldsymbol{n}=\boldsymbol{n}^{+}=-\boldsymbol{n}^{-}$.
We combine the similar items

\begin{equation}
	\begin{split}\int_{\Omega^{+}}a^{+}(\nabla u)\cdot(\nabla\delta u)d\Omega+\int_{\Omega^{-}}a^{-}(\nabla u)\cdot(\nabla\delta u)d\Omega\\
		-\int_{\varGamma^{t+}}\bar{t}\delta udx-\int_{\varGamma^{t-}}\bar{t}\delta ud\varGamma+\int_{\varGamma^{inter}}(u^{+}-u^{-})\delta ud\varGamma & =0.
	\end{split}
	\label{eq:energy delta}
\end{equation}
The  \Cref{eq:energy delta} is equal to the energy stationary, so
we can get the energy form with the interface,
\begin{equation}
	\begin{split}J & =\frac{1}{2}\int_{\Omega^{+}}a^{+}(\nabla u)\cdot(\nabla u)d\Omega+\frac{1}{2}\int_{\Omega^{-}}a^{-}(\nabla u)\cdot(\nabla u)d\Omega\\
		& -\int_{\varGamma^{t+}}\bar{t}ud\varGamma-\int_{\varGamma^{t-}}\bar{t}ud\varGamma+\beta\int_{\varGamma^{inter}}(u^{+}-u^{-})^{2}d\varGamma
	\end{split}
\end{equation}
where $\beta$ is the hyperparameter of the interface. Further,  the above energy form is a convex function, and the exact
solution is the stationary point, i.e., the extreme value problem. Although the energy form is convex in terms of
the whole function space, it is commonly not convex in the neural
network function space. So we can get the solution by optimizing
J to a minimum. We can assign the different neural network in each subdomain.
The backpropagation of CENN is discussed in \ref{sec:Appendix-B.-Back propogation of CENN}.
CENN penalty only has one unknown penalty on the interface to ensure the
continuity of the interesting field. However,  CPINN not only has
the $\boldsymbol{u}$ continuity condition but also has the derivative of $\boldsymbol{u}$
continuity condition. The additional derivative of CPINN will decrease
the accuracy and efficiency. If the order of the PDE increase,  the
additional term about interface will increase in CPINN more than CENN,
e.g., if the order of PDE is $2m$ order derivative,  it is easy to
prove that CPINN will have m more interface penalty factors than CENN.
Although there are currently some  theoretical guidances for the selection of hyperparameters \citet{NTK_PINN,NTK_to_get_hyperparameter_of_PINN,ill_gradient}, 
there are still many challenges to determine the hyperparameters, i.e., the different components of the loss function.  Compared with the CPINN, the penalty of CENN is greatly reduced,
 but there is still a penalty of the interface,  which is a hyperparameter.
Here we use a heuristic construction 
\begin{equation}
	\beta=-c\cdot ln(tanh(\frac{N_{interface}}{N_{domain}}))\label{eq:interface_hyper}
\end{equation}
 where  $c$ is
a scale factor, we recommend 1e3. $N_{interface}$ and $N_{domain}$
are the number of points at the interface and the domain respectively.
The \Cref{eq:interface_hyper} uses the concept of information entropy,  we assume
that $p=tanh(N_{interface}/N_{domain})$ is the probability of the
certainty of the interface.  Obviously,  the more interface points,
 the greater certainty of probability about the interface. 
The training is essentially  multi-task learning,  therefore it is a game between
the interface loss and the energy functional.  If there are more interface
points,  less attention has been paid to energy functional training.  We need to adjust the hyperparameter $\beta$ to balance the
interface loss and energy loss.  Finally,  $-ln(p)$ is the information
entropy to evaluate the value of the interface information, i.e.,
it indicates that the interface has more abundant information if there
are more interface points. The above heuristic construction is used
to construct the penalty of the interface.  After many numerical experiments,
 we found that training is successful by this construction of the
hyperparameter. 

\begin{rem}
CENN is suitable for solving the problem of one point
coordinate containing multiple values of interesting field, such as the crack
problem,  which will be mentioned in detail in \Cref{sec:crack problem}.
CENN is also suitable for solving the problem of weak continuity problem,
i.e.,  the original function is continuous but the
derivative function is discontinuous,  such as heterogeneous problem,
which will be mentioned in detail in \Cref{sec:non homogeneous boundary problem} and
\Cref{sec:hpyerelasticity material}. 
\end{rem}
\begin{rem}
	DEM is a special form of the weighted residual. The test functuion
	of DEM  and traditional PINN is $\partial u(\boldsymbol{x};\boldsymbol{\theta})/\partial \theta_{j}$ and  $2\triangle\partial u(\boldsymbol{x};\boldsymbol{\theta})/\partial\theta_{j}$ respectively.
\end{rem}
For the sake of simplicity,  we consider the Poisson equation
with the whole essential boundary 
\begin{equation}
	\triangle u=f.
\end{equation}
The energy form of the strong form is 

\begin{align}
	\mathcal{L} & =\sum_{i=1}^{n^{v}}w_{i}w_{\varepsilon}(u_{1}(\boldsymbol{x}_{i};\boldsymbol{\theta}))-\sum_{i=1}^{n^{v}}w_{i}f_{i}u(\boldsymbol{x}_{i};\boldsymbol{\theta})
\end{align}
where $w_{i}$ is the weight of the attribution points $\boldsymbol{x}_{i}$,  especially
$w_{i}=V/n$ if uniform random Monte Carlo method is adopted.  $w_{\varepsilon}=\frac{1}{2}(\nabla u)\cdot(\nabla u)$,
$u$ is a scalar interesting variable.  The first-order variation w.r.t.  Loss is 
\begin{align}
	\delta \mathcal{L} & =\sum_{i=1}^{n^{v}}\frac{V}{n}(\nabla u_{i})\cdot(\nabla\delta u_{i})-\sum_{i=1}^{n^{v}}\frac{V}{n}f_{i}\delta u_{i}
\end{align}
The first term on the RHS is integrated by parts.  We can get
\begin{equation}
	\begin{split}\delta \mathcal{L} & =\sum_{i=1}^{n^{v}}\frac{V}{n}(\triangle u_{i})(\delta u_{i})-\sum_{i=1}^{n^{v}}\frac{V}{n}f_{i}\delta u\\
		& =\sum_{i=1}^{n^{v}}\frac{V}{n}(\triangle u_{i}-f_{i})(\delta u_{i}).
	\end{split}
\end{equation}
Given the trial function $u$ is an admissible function,  the boundary
part generated by integrating by parts with $\delta u$ is vanish.
Considering the interesting field $u$ is the function of neural
network parameter $\theta$,  the variation form can be rewritten
\begin{equation}
	\begin{split}\delta \mathcal{L} & =\sum_{i=1}^{n^{v}}\frac{V}{n}(\triangle u_{i}-f_{i})(\sum_{j=1}^{n^{\theta}}\frac{\partial u_{i}}{\partial\theta_{j}}\delta\theta_{j})\\
		& =\sum_{j=1}^{n^{\theta}}\sum_{i=1}^{n^{v}}\frac{V}{n}(\triangle u_{i}-f_{i})(\frac{\partial u_{i}}{\partial\theta_{j}}\delta\theta_{j}).
	\end{split}
\end{equation}
The stationary point (the first-order variation is zero)
is equal to the strong form with the test function $\partial u(x;\theta)/\partial\theta_{j}$,
j=1, 2, . . . $n^{\theta}$,  i. e. , 
\begin{equation}
	\begin{split}\mathcal{L}_{1} & =\sum_{i=1}^{n^{v}}\frac{V}{n}(\triangle u_{i}-f_{i})(\frac{\partial u_{i}}{\partial\theta_{1}})\\
		\mathcal{L}_{2} & =\sum_{i=1}^{n^{v}}\frac{V}{n}(\triangle u_{i}-f_{i})(\frac{\partial u_{i}}{\partial\theta_{2}})\\
		\vdots\\
		\mathcal{L}_{n^{\theta}} & =\sum_{i=1}^{n^{v}}\frac{V}{n}(\triangle u_{i}-f_{i})(\frac{\partial u_{i}}{\partial\theta_{n^{\theta}}})
	\end{split}
\end{equation}
DEM can be thought of as jointly optimizing all the above losses so
that all losses converge to zero.

On the other hand, we analyze the traditional PINN, i.e., the least
square method
\begin{equation}
	\mathcal{L}=\sum_{i=1}^{n^{v}}\frac{V}{n}(\triangle u_{i}-f_{i})^{2}.
\end{equation}
The first order variation w. r. t. PINN Loss is 
\begin{align}
	\delta\mathcal{L} & =\sum_{i=1}^{n^{v}}\frac{V}{n}2(\triangle u_{i}-f_{i})(\delta\triangle u_{i}).
\end{align}
Considering the interesting field $u$ is the function of neural network
parameter $\theta$, we can obtain
\begin{align}
	\delta\mathcal{L} & =\sum_{j=1}^{n^{\theta}}\sum_{i=1}^{n^{v}}\frac{V}{n}2(\triangle u_{i}-f_{i})(\frac{\partial\triangle u_{i}}{\partial\theta j}\delta\theta_{j}).
\end{align}
The traditional PINN is equal to the strong form with the test function
$2\partial\triangle u/\partial\theta_{j}$ , j=1, 2, . . . $n^{\theta}$,
i. e. , 
\begin{equation}
	\begin{split}\mathcal{L}_{1} & =\sum_{i=1}^{n^{v}}\frac{V}{n}(\triangle u_{i}-f_{i})(2\frac{\partial\triangle u_{i}}{\partial\theta_{1}})\\
		\mathcal{L}_{2} & =\sum_{i=1}^{n^{v}}\frac{V}{n}(\triangle u_{i}-f_{i})(2\frac{\partial\triangle u_{i}}{\partial\theta_{2}})\\
		\vdots\\
		\mathcal{L}_{n^{\theta}} & =\sum_{i=1}^{n^{v}}\frac{V}{n}(\triangle u_{i}-f_{i})(2\frac{\partial\triangle u_{i}}{\partial\theta_{n^{\theta}}})
	\end{split}
\end{equation}
Traditional PINN can be thought of as jointly optimizing all the above
losses so that all losses converge to zero.

So the test functuion of DEM is $\partial u(\boldsymbol{x};\boldsymbol{\theta})/\partial\theta_{j}$,
the test function of traditional PINN is $2\partial\triangle u(\boldsymbol{x};\boldsymbol{\theta})/\partial\theta_{j}.$
In fact, both of
these are actually special cases of the weighted residual method.

If the second order coordinate derivative of the interesting field
is not close to zero, though first order derivative close to zero,
e.g., the minimum value problem, the strong form may be better than
energy form because test function is not zero compared to energy form. 

\section{Result\label{sec:Result}}

\subsection{Crack\label{sec:crack problem}}

\begin{figure}
\begin{centering}
\subfloat{\centering{}\includegraphics{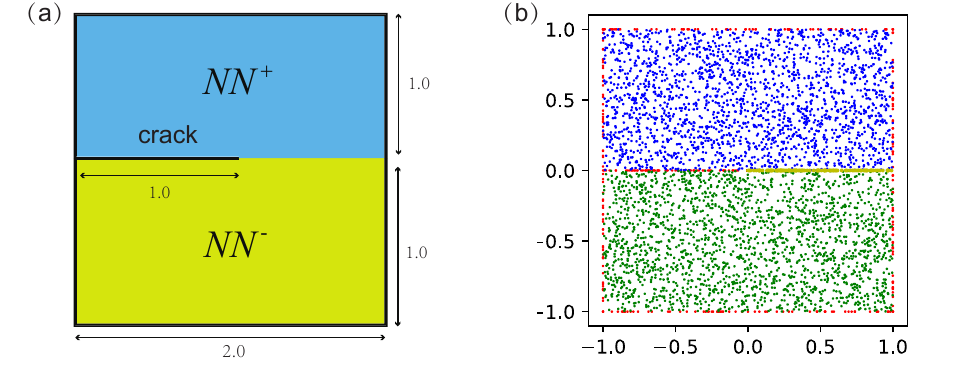}}
\par\end{centering}
\caption{ The schematic diagram of \uppercase\expandafter{\romannumeral3} mode
crack: (a) NN$^{+}$ represents the neural network in the upper area.
NN$^{-}$represents the different neural network in the lower area.
the geometric size is a square with side length 2, the length of the crack
is 1, spanning in the middle of the structure , i.e., x<0, y=0. (b)
 Illustration of different sampling sizes and strategies, the blue points
are the internal training points of the neural network $NN^{+}$,
the green points are the internal training points of the neural network
$NN^{-}$. The red points are the training points of the essential
boundary conditions, and the yellow points are the interface points
of the different neural network. The total number of internal points
in the upper and lower regions is 4096, the number of essential boundaries
condition points is 256, and the number of interface points is 1000.
All points are randomly distributed unless otherwise stated.
\label{fig:crack schematic}}
\end{figure}

This subsection shows the proposed method CENN can tackle the strong discontinuity and singularity problem. The governing equation of \uppercase\expandafter{\romannumeral3}
mode crack is given by
\begin{align}
	\triangle(u(\boldsymbol{x}))=0 & ,\boldsymbol{x}\in\Omega\label{eq:laplace equation}
\end{align}
where $\Omega=(-1,1)\vartimes(-1,1)$. The analytical displacement
solution is $u=r^{\frac{1}{2}}sin(\frac{1}{2}\theta)$,  $r$ and
$\theta$ are the radius and the angle of polar coordinates respectively, and the
$\theta$ range is $[-\pi,+\pi]$.  We construct the boundary conditions
through the analytical solution. It is worth noting that the solution
to this problem suffers from the well-known ``corner singularity'',
 which means that the derivative is one-half singularity at the center
point(0, 0) \citet{finite_element_book}. Given that the displacement
solution of this problem is discontinuous at the crack (x<0,  y=0),
 the same coordinate point at the crack corresponds to multiple displacement
values $\pm r^{\frac{1}{2}}$, so it is necessary to use two neural
networks to approximate the displacement field. We divide the area
into upper and lower parts,  and use different neural networks to
approximate them,  as shown in \Cref{fig:crack schematic}a. We consider the commonly used data-driven and CPINN (PINN strong form with subdomains) methods for comparison with CENN. Data-driven is using neural nerwork to fit analytical or reliable numerical solution calculated in advance, such as FEM. The input is the coordinates; The output is the field of interest. The loss function (usually MSE) is the difference between the neural network prediction and the field of interest that has been calculated in advance.
These three  methods have the same point allocation way,  as shown in \Cref{fig:crack schematic}b.  Training points are redistributed every 100 epochs in all
three methods. The energy principle needs to convert the  strong
form into a variational form.  The variational form of the \Cref{eq:laplace equation}
is
\begin{equation}
J(u)=\int_{\Omega}\frac{1}{2}(\nabla u)\cdot(\nabla u)d\Omega.\label{eq:crack functional}
\end{equation}
Using the energy principle,  the minimum value of the above formula
is equivalent to the strong form of the solution,  but the energy
principle needs to construct a possible displacement field, i.e.,
the admissible function. It is necessary to find an exact displacement
in the possible displacement field,  so that the functional of \Cref{eq:crack functional}
takes the minimum value,  and the displacement is considered to be the
optimal solution in the sense of optimal energy error, i.e., $\underset{u^{pred}}{argmin}\int \varPsi(u^{error})d\Omega$,
where $u^{error}=u^{*}-u^{pred}$ and $\varPsi$ is energy density; $u^{*}$ and $u^{pred}$ is the exact value and prediction of
the problem \citet{finite_element_book}.

\begin{figure}
\begin{centering}
\includegraphics[scale=0.85]{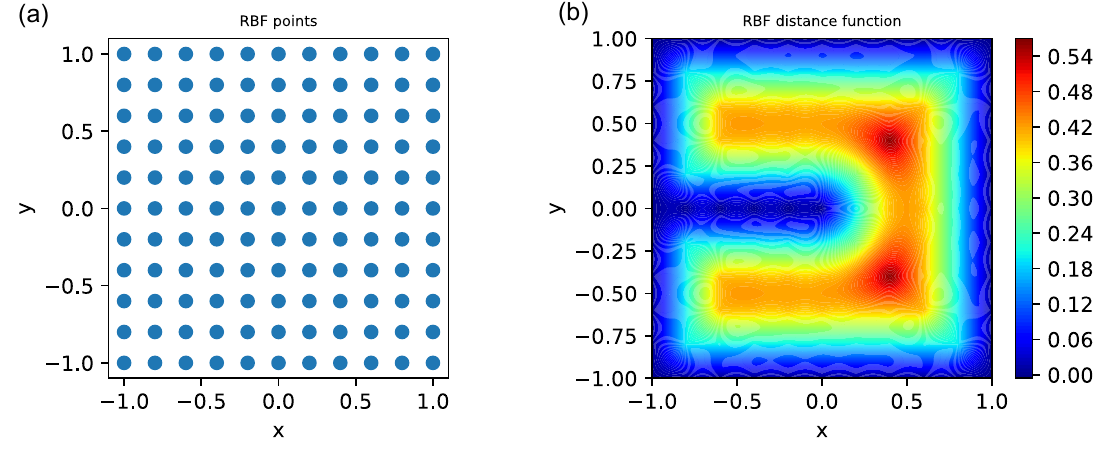}
\par\end{centering}
\caption{The RBF  of \uppercase\expandafter{\romannumeral3} mode
	crack. (a) The RBF distance network of \uppercase\expandafter{\romannumeral3}
mode of the crack: RBF points allocation map. (b) RBF distance
network prediction contour. \label{fig:RBF}} 
\end{figure}

\begin{figure}
\begin{centering}
\includegraphics[scale=0.85]{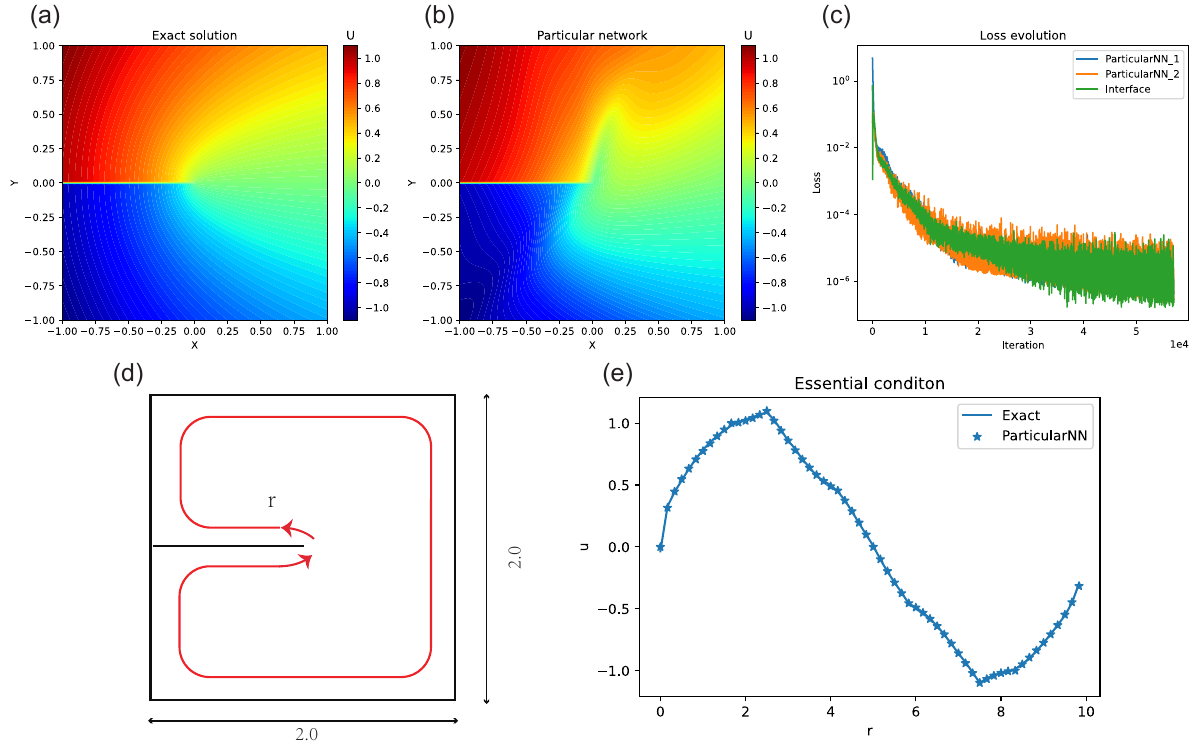}
\par\end{centering}
\caption{The particular neural network  of \uppercase\expandafter{\romannumeral3} mode crack. (a) The contour of the exact solution of \uppercase\expandafter{\romannumeral3}
mode crack. (b) The contour of the particular prediction. (c) Evolution of the MSE loss for  training particular
network, particularNN\_1 and 2 represent the neural network satisfying the essential boundary condition of up and down region respectively. The number and the arrangement of points
in both neural networks are same. Points are redistributed every 10
epoch. (d) The evaluation trace of the particular solution,
starting from the center point and going around the boundary clockwise. (e) The comparison of the particular network and the exact solution
on the trace of the boundary clockwisely\label{fig:particular neural network of crack}}
\end{figure}

Therefore,  compared with CPINN,  CENN requires additional
construction of possible displacement fields that meet the essential
boundary conditions.  The RBF distance network and the particular
network that meets the essential boundary conditions need to be constructed
in advance. The number of RBF allocation points and the arrangement
of points need to be determined according to the geometric shape of
the essential boundary conditions. In this example,  there are 121 uniform points in the domain as shown in \Cref{fig:RBF}a.
RBF distance network prediction is shown in  \Cref{fig:RBF}b.
The average error of the RBF network is about 0.76 \%.  The accuracy
can be improved by increasing the RBF center point including boundary
points and the domain points, but it will increase the amount of calculation.
 After obtaining the RBF distance network, we train two particular
networks in the upper and lower region to fit the essential boundary
conditions.  The hidden layer of the particular network has 3 layers,
 each layer has 10 neurons more shallow than the general network, the
activation function is tanh in all the hidden layers, and an identical function
is used in the output layer,  the optimizer is Adam \citet{ADAM}, the 
learning rate is 5e-4.  It is worth noting that the neural network
structure of the particular solution network does not need to be too
complicated thanks to the powerful fitting ability of the neural network.
 A simple neural network structure can fit the essential boundary
conditions well.  It is worth noting that we must add $\mathcal{L}_{interface}=||u^{+}-u^{-}||$
at the interface to satisfy the continuity,  and the convergence condition
of the particular network training is that the MSE is less than 1e-6.
\Cref{fig:particular neural network of crack}a and b show that the output of the particular network is close to the pattern
of the exact solution because the boundary conditions are all essential
boundary conditions enclosing the interesting domain.  As shown in
\Cref{fig:particular neural network of crack}c,  the
network structure can fit the boundary value well as the number of
network training rounds increases. The loss function of the boundary
drops a little faster than the loss value of the interface because
there is a specific label at the boundary.  The convergence
speed of the boundary loss of different particular networks is almost
the same,  which shows that the particular network can learn the x-axis
center symmetry of the exact solution very well. It does not matter that 
the particular network and the exact solution are quite different
in the domain without essential boundaries.  Because the function of
the particular network is to precisely satisfy the given value at
the essential boundary conditions. To further explore the
performance of the particular network,  we make a clockwise circle
around the essential boundary as shown in  \Cref{fig:particular neural network of crack}d. \Cref{fig:particular neural network of crack}e shows the comparison of the particular network and the
exact solution at the essential boundary.  We can find that the particular
network can fit the essential boundary conditions very well,  which
 is the result of the powerful fitting ability of the neural network. 

After constructing the possible displacement field, we compare data-driven,
 CPINN and CENN.  The number of allocation points of these three ways
is the same.  The network structure and optimization scheme are also
the same.  It is worth noting that CPINN does not use the construction
of the admissible function, but the other two methods (data-driven
and CENN) do that. The generalized network
has 4 layers,  each layer has 20 neurons,  and the learning rate is
1e-3.  We consider the admissible function

\begin{equation}
	\boldsymbol{u}(\boldsymbol{x})=\boldsymbol{u}_{p}(\boldsymbol{x};\boldsymbol{\theta}_{p})+RBF(\boldsymbol{x})\cdot \boldsymbol{u}_{g}(\boldsymbol{x};\boldsymbol{\theta}_{g})
\end{equation}
where $\boldsymbol{u}_{p}$ is the particular network, $\boldsymbol{\theta}_{p}$ is the parameter
of the particular network,  RBF is the radial basis function,  $\boldsymbol{u}_{g}$
is the generalized network,  $\boldsymbol{\theta}_{g}$ is the parameter of the generalized
network. In the training process,  we freeze the parameters of the
particular network and RBF,  and only train the generalized neural
network through the gradient descent method, i.e., Adam.  It is worth
noting that RBF is equal to zero at the essential boundary, so the
role of the generalized network is eliminated, which results that
the particular network only plays games. If the training of the RBF
and particular network is successful, then the admissible function accurately
meets the given displacement value at the essential boundary,  which
is the heart of the admissible function.  On the other hand,  the
CPINN subdomains are the same as CENN, i.e., the upper and lower neural
networks. The boundary conditions and the interface conditions are
imposed by the penalty method, noting that the interface conditions
include not only the continuous condition of displacement but also
the continuous condition of displacement derivative, i.e.,
\begin{equation}
	\begin{split}\boldsymbol{n}\cdot(\nabla u^{+}-\nabla u^{-})=0, & \boldsymbol{x}\subseteq\varGamma^{inter}\\
		u^{+}(\boldsymbol{x})=u^{-}(\boldsymbol{x}), & \boldsymbol{x}\subseteq\varGamma^{inter}.
	\end{split}
\end{equation}
 Compared to CENN,  one more loss is added,  and the additional loss
function is the derivative form,  which will reduce the accuracy and
efficiency. The loss function of CPINN is
\begin{equation}
\begin{split}\mathcal{L} & =\lambda_{1}\int_{\Omega^{+}}|\triangle(u^{+})|^{2}d\Omega+\lambda_{2}\int_{\Omega^{-}}|\triangle(u^{-})|^{2}d\Omega+\lambda_{3}\int_{\varGamma^{+}}|u^{+}-\bar{u}|^{2}d\varGamma+\lambda_{4}\int_{\varGamma^{-}}|u^{-}-\bar{u}|^{2}d\varGamma\\
 & +\lambda_{5}\int_{interface}|u^{-}-u^{+}|^{2}d\varGamma+\lambda_{6}\int_{interface}|\boldsymbol{n}\cdot(\nabla u^{+}-\nabla u^{-})|^{2}d\varGamma
\end{split}
\end{equation}

\begin{figure}
\begin{centering}
\includegraphics[scale=0.3636]{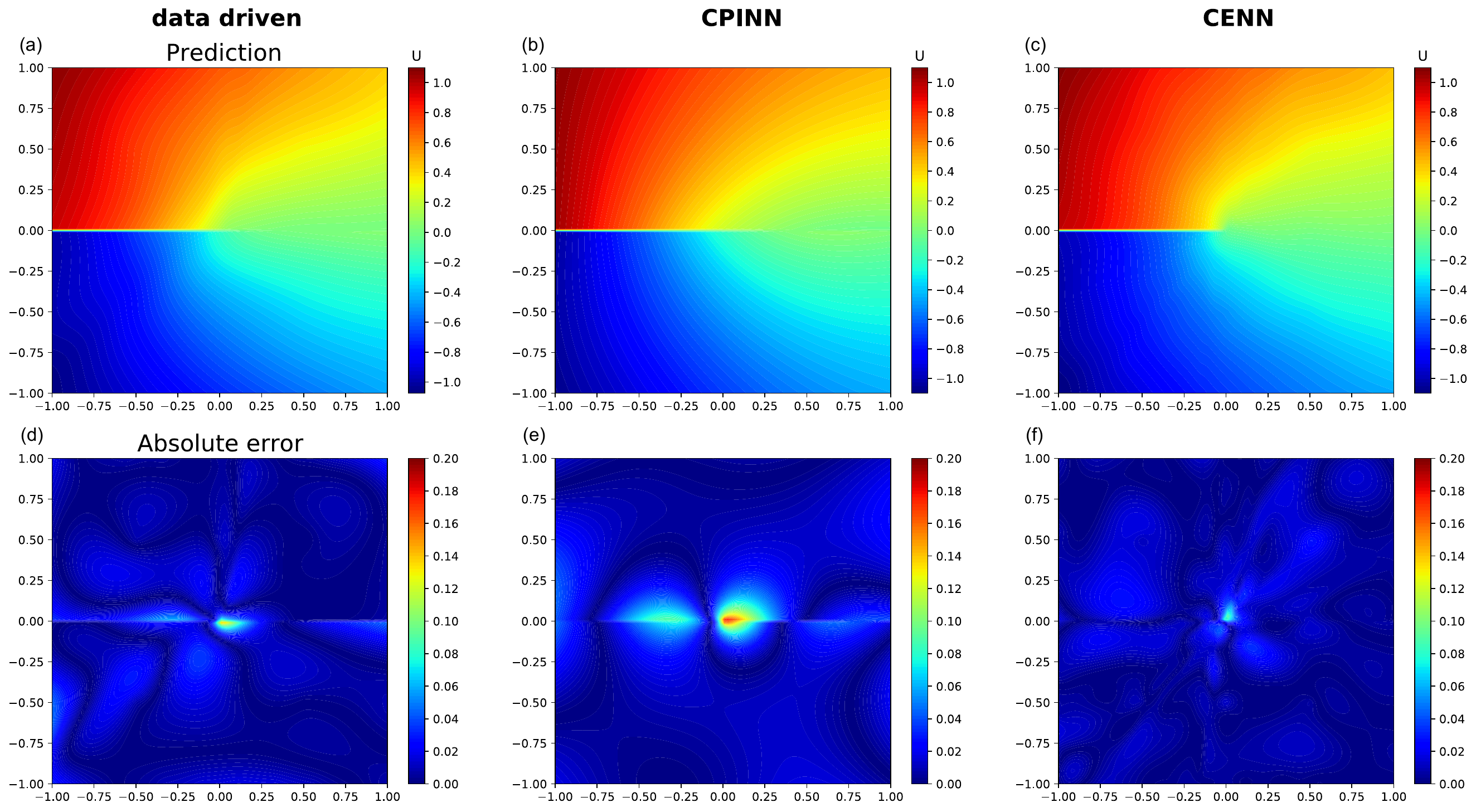}
\par\end{centering}
\caption{  Purely date-driven model with mean square error (MSE) of displacement field: prediction (a) and absolute error (d). CPINN model with  MSE loss function:  prediction (b) and absolute error (e). CENN model with energy functional loss and $\lambda_{interface}=2700$: prediction (c) and absolute error (f).
\label{fig:Comparison_predition_error_crack}}
\end{figure}

\begin{figure}
\begin{centering}
\includegraphics[scale=0.2936]{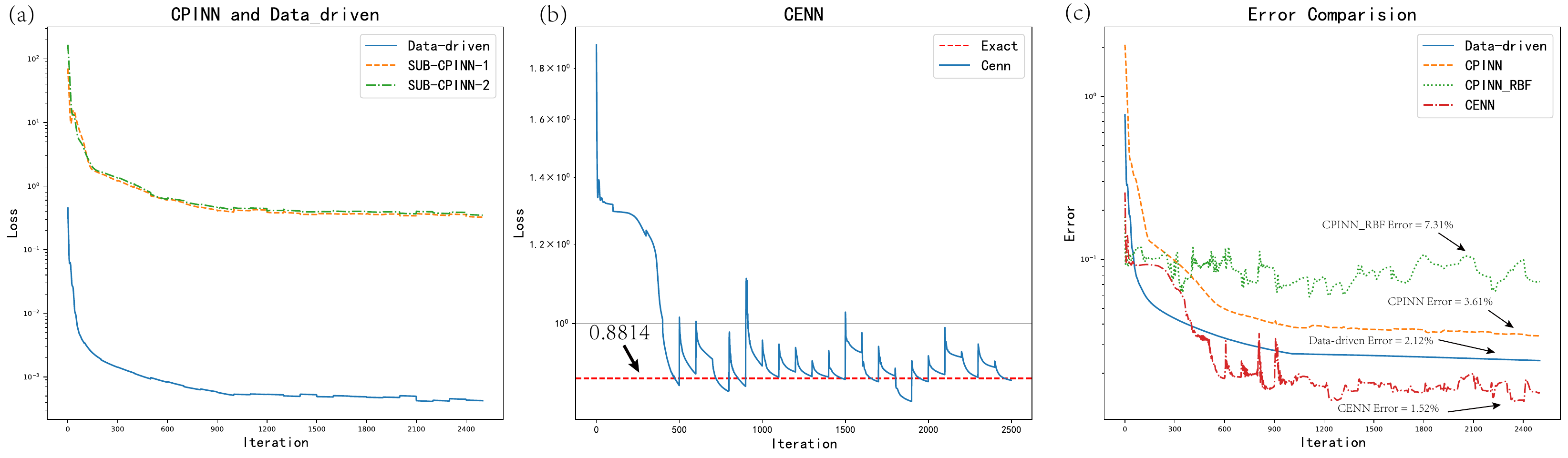}
\par\end{centering}
\caption{Comparison the loss and error of data-driven, CPINN and CENN in \uppercase\expandafter{\romannumeral3}
mode crack: (a) Comparison of the  CPINN loss value  and data-driven loss value. SUB-CPINN-1 and 2 are CPINN subdomains NN in upper region and lower region of crack respectively. (b) Comparison of  CENN energy functional  and  exact functional
integration. (c) Comparison of the overall relative error $L_2$ of data driven, CPINN, CPINN-RBF and CENN.\label{fig:Comparison_position_crack}}
\end{figure}

\begin{figure}
\begin{centering}
\includegraphics[scale=0.5]{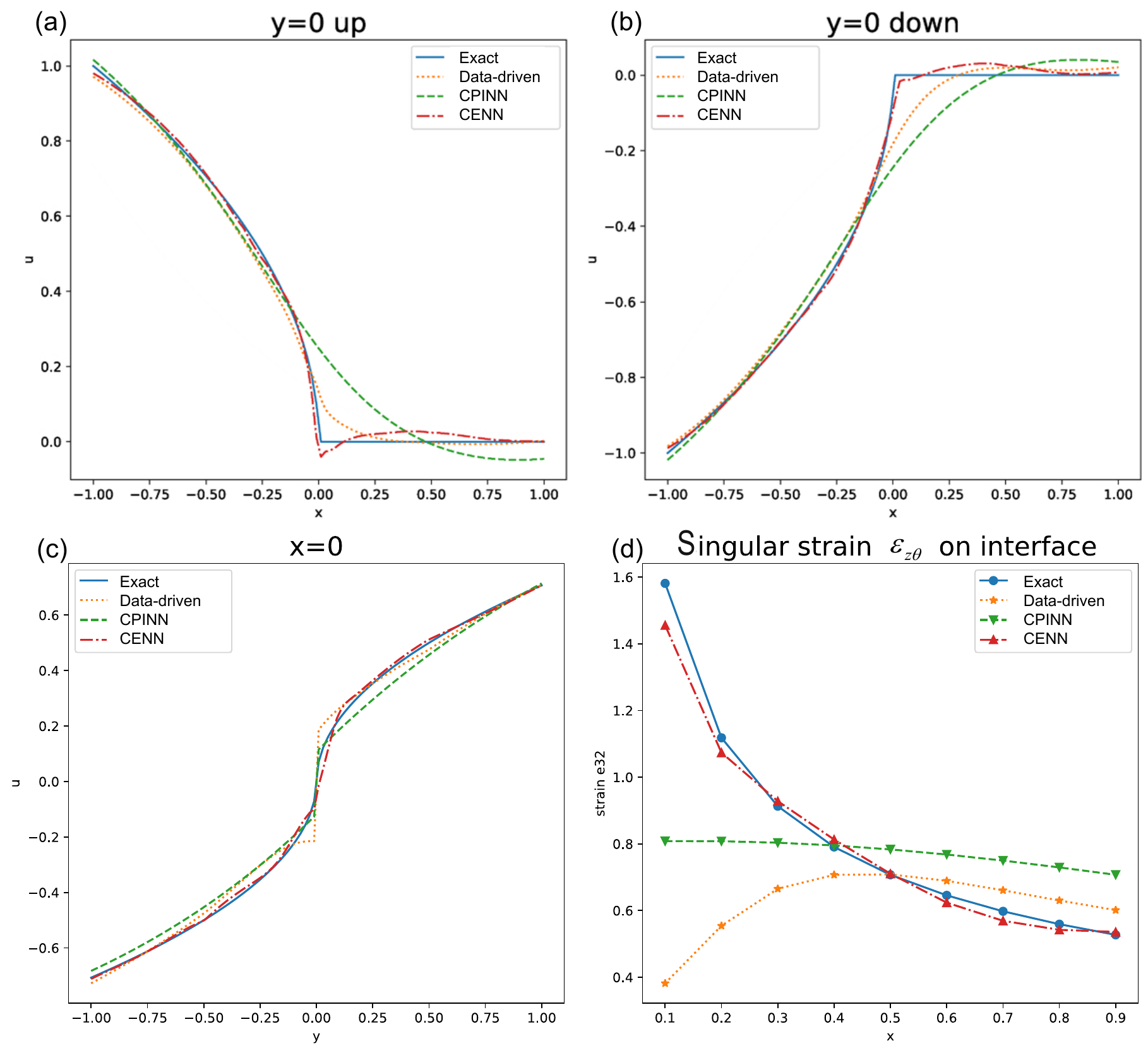}
\par\end{centering}
\caption{Different location comparison among data-driven, CPINN, and CENN: (a) Comparison
of the displacement solution on the horizontal line y=0 up. (b) Comparison of the displacement solution on the horizontal line
y=0 down. (c) Comparison of displacement solution at x=0. (d) Comparison of singular displacement derivative ${\partial u}/{\partial y}$ ($\varepsilon_{z\Theta}$)
at x>0, y=0.\label{fig:Comparison of different locations of cracks}}
\end{figure}

It is not difficult to find that there are many hyperparameters in
CPINN (more subdomains will further increase hyperparameters).  Through
tuning repeatedly,  we have selected the best set of hyperparameters,
 $\lambda_{1}=\lambda_{2}=1,\lambda_{3}=\lambda_{4}=50,\lambda_{5}=\lambda_{6}=10$
(maybe there are better hyperparameters,  because there are 6 hyperparameters,
 the workload of tuning parameters is very large).  \citet{NTK_to_get_hyperparameter_of_PINN}
shows that it is possible to use NTK (Neural tangent kernel) theory to automatically adjust
hyperparameters for PINN,
 but there is no method for hyperparameter adjustments in the form
of subdomain PINN, i.e., CPINN. \Cref{fig:Comparison_predition_error_crack}a and e show the data-driven prediction results and absolute error,
\Cref{fig:Comparison_position_crack}c shows that the overall relative error of the data-driven is
2.12 \%. Due to the randomness of the initial parameters of the neural
network,  the initial parameters are initialized with Xavier \citet{Xavier},
i.e. a method to initialize parameters. If not specified,  all the
results below are the average of 5 times.  Since data-driven does
not consider PDE,  the loss function does not include auto differential
\citet{automatic_differential},  which has higher accuracy and computational
efficiency.  We found the absolute error near the center (x=0,  y=0)
is the largest as shown in \Cref{fig:Comparison_predition_error_crack}d. \Cref{fig:Comparison of different locations of cracks}a and b show that the exact solution is not smooth at the center.
This is the reason for the large error at the center point.  In
addition,  there were relatively sharp fluctuations at the beginning
loss evolution, and then gradually converged. \Cref{fig:Comparison_predition_error_crack}b and e shows that the prediction results of CPINN are sharp relatively,
 and the errors are also mainly concentrated in the center point.
 Because the derivative of analytical solution with respect to $\theta$
is singular at the center point,  causing the loss at the center point
to fluctuate sharply.  Considering that CPINN and data-driven are
both MSE errors, the optimal loss values both are 0,  so we compare
the loss functions of the two together as shown in \Cref{fig:Comparison_position_crack}a.  In addition,  CENN is based on the energy method,  and the
loss function is an energy functional,  so the optimal loss value
is not 0. From \Cref{fig:Comparison_position_crack}a,  it can be seen that the CPINN loss function is higher than
data-driven.  The CPINN loss function is unstable compared to the
data-driven,  which is due to the high-order derivatives involved.
 In addition,  the decreasing trend of the loss function of the different
networks is similar in CPINN.  The reason is 
the x-axis symmetry of the crack problem.  To further analyze the performance of these
three methods,  we compare the relative error $L_{2}=\sqrt{\Sigma_{i=1}^{N_{pred}}error^{2}(\boldsymbol{x}_{i})}/\sqrt{\Sigma_{i=1}^{N_{pred}}exact^{2}(\boldsymbol{x}_{i})}$
 in all three methods. We include the comparison of CENN and CPINN-RBF (CPINN: Boundary conditions are satisfied by penalty method, CPINN-RBF: Boundary conditions are satisfied by constructing the admissible function with RBF in advance). \Cref{fig:Comparison_position_crack}c shows that the overall relative error of CPINN is 3.61 \%, and the overall relative error of CPINN-RBF is 7.31 \%. 
 This is one of the reasons that the weight of the different loss
functions is not adjusted to the optimal value. In addition,  CPINN
involves higher-order derivatives,  so the efficiency and accuracy
are not as good as CENN. \Cref{fig:Comparison_predition_error_crack}c and f show that the prediction result of CENN is smooth,  not as
sharp as CPINN. \Cref{fig:Comparison_position_crack}c shows that the overall relative error of CENN is 1.52 \% and
the overall result of CENN is more precise than CPINN and even better than data-driven
accuracy.  Due to the low derivative order of the loss function,  the
efficiency is higher than CPINN.  The error distribution is more uniform.
The main error is at the center point (x=0,  y=0),  which is caused
by the form of the admissible function. This will be analyzed in
the discussion part in \Cref{sec:Discussion} in detail.  \Cref{fig:Comparison_position_crack}b shows that CENN loss function converges to the
exact value of the functional

\[
J=\int\frac{1}{2}\text{[}(\frac{\partial u}{\partial r})^{2}+(\frac{\partial u}{r\partial\theta})^{2}]rdrd\theta=\int[\frac{1}{8\sqrt{x^{2}+y^{2}}}]dxdy.
\]
The analytic integration of the energy functional is about
0. 8814.  According to the principle of minimum potential energy,
 the exact solution is the minimum value of the energy functional,
i.e., 0.8814 in this problem. The reason for a small amount of fluctuation
in the later stage of the training is due to the learning rate and the numerical integration accuracy, as shown in \Cref{fig:Comparison_position_crack}b.  This
can be eliminated by reducing the learning rate and adopting a more precise numerical integration scheme.  In addition,  the
loss function is sometimes slightly lower than the exact value due
to numerical integration and discrete errors, which can be reduced
by increasing random points.  \Cref{fig:Comparison_position_crack}c shows that CENN outperforms CPINN obviously,  and even surpasses
data-driven that does not involve coordinate derivatives.  Finally,  we can
find that the trends of the loss function and relative error of the
above three methods are the same. 

Next,  we compare the discontinuous displacement solutions and singular
strains that we are more concerned about.  \Cref{fig:Comparison of different locations of cracks}
shows the comparison of solution of data-driven, CPINN and CENN at
various x and y position.  We observe the displacement of y=0,  the
up and down displacement.  This is because at (x<0, y=0),  the displacement
solution is multiple values of the same coordinate,  so there are
two displacement solutions as shown in \Cref{fig:Comparison of different locations of cracks}a,b.  We take into account the strain
\begin{equation}
\varepsilon_{z\Theta}|_{interface}=\frac{1}{r}\frac{\partial u}{\partial\theta}=\frac{1}{2\sqrt{r}}cos(\frac{\theta}{2})|_{\theta=0}=\frac{1}{2\sqrt{r}},
\end{equation}
the y-direction derivative of the displacement is singular at the
center (x=0,  y=0),  so we research the singular strain at the interface,
 and explore whether the neural network is capable of fitting the
singularity problem.  \Cref{fig:Comparison of different locations of cracks}a and b shows that the prediction of CENN and data-driven is close to
the displacement solution,  and both have higher accuracy than CPINN.
 The error between CENN and data-driven is large at x=0 caused by the RBF distance network of the admissible function,  which
will be further analyzed in \Cref{sec:Discussion}.  \Cref{fig:Comparison of different locations of cracks}c shows the comparison of x=0.  The accuracy of CENN at
x=0 is higher than that of data-driven as the same as y=0, and the
error near the center point is higher.  It is also the reason that
the RBF distance network of the admissible function causes. \Cref{fig:Comparison of different locations of cracks}d
shows the singular strain.  The performance of CENN in singular strain
is better than data-driven and CPINN.  This is because loss
function of CENN is constructed based on the physical minimum energy
theory whose interpretability is stronger.  The strain will cause
the loss function to change in CENN,  which will make the CENN strain
converge to the exact solution.  In numerical methods such as finite
element,  the singular point of fracture mechanics often requires a 
special quarter-node element \citet{finite_element_book},  but this
method does not have a special treatment,  which comes from the tremendous
function space of the neural network itself.

\subsection{Non homogeneous problem with  complex boundary \label{sec:non homogeneous boundary problem}}

\begin{figure}
\begin{centering}
\includegraphics[scale=0.65]{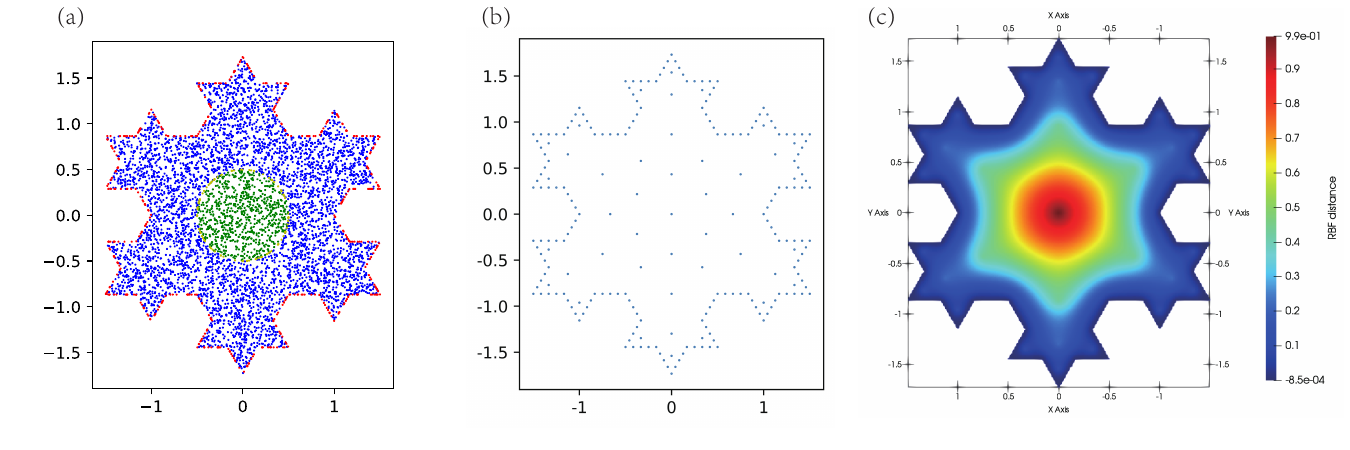}
\par\end{centering}
\caption{Distance network with complex boundaries of inhomogeneity: (a) Illustration
of different sampling size and points allocation strategy, the green
points are the PDEs parameter $a=$$1/15$($r<0.5$) training points, the blue points
are the $a=1$($r>0.5$) training points, the red points are the training
points of the essential boundary conditions, and the yellow points
are the interface points ($r=0.5$) of the different domain, the total
number of domain points is 10000, the number of essential boundary
condition points is 4800, and the number of interface points is 1000.
All points are randomly distributed. (b) The collocation point
diagram of RBF network, uniform distribution, 217 points on the boundary, 48 points in the
domain. (c) Prediction of RBF distance network.\label{fig:distance koch}}
\end{figure}

In this section, we investigate the problem with the complex boundary to shows the proposed method CENN can tackle the  complex boundary and heterogeneous problem. Complex boundary problem can not be solved well in a traditional method such as FEM \citet{complex_PINN_a_method_to_construct_admissible_function}.
 So we consider the Koch snowflake (complex boundary problem)  as our boundary. We also solve
the non-homogeneous problem,  which is widespread in physics such
as composites material,  as shown in \Cref{fig:distance koch}a,  the fractal level L=2.   The govern equation of the problem is \citet{PINN_energy_form_to_solve_C0_without_subdomains}

\begin{equation}
	\begin{cases}
		-a(\boldsymbol{x})\triangle u(\boldsymbol{x})=f(\boldsymbol{x}) & \boldsymbol{x}\in\Omega\\
	u(\boldsymbol{x})=g(\boldsymbol{x}) & \boldsymbol{x}\in\partial\Omega,
	\end{cases}
\end{equation}
where $\alpha(x)=\alpha_{i},x\in V_{i}$,  $\alpha$ for different
regions is a different constant.  For the sake of simplicity,  we consider
\begin{equation}
\begin{cases}
a_{1}=\frac{1}{15} & r<r_{0}\\
a_{2}=1 & r\geq r_{0}.
\end{cases}
\end{equation}
 We adopt the exact solution
\begin{equation}
u(r,\theta)=\begin{cases}
\frac{r^{4}}{\alpha_{1}} & r<r_{0}\\
\frac{r^{4}}{\alpha_{2}}+r_{0}^{4}(\frac{1}{\alpha_{1}}-\frac{1}{\alpha_{2}}) & r\geq r_{0}.
\end{cases}
\end{equation}
 The energy form of the problem is 
\begin{equation}
J=\int_{\Omega_{1}}\frac{1}{2}\alpha_{i}(\nabla u)\cdot(\nabla u)d\Omega+\int_{\Omega_{2}}\frac{1}{2}\alpha_{i}(\nabla u)\cdot(\nabla u)d\Omega-\int_{\Omega}f ud\Omega.
\end{equation}
where $\Omega_{1}$ is the domain $r<r_{0}$ and $\Omega_{2}$ is
the domain $r\geq r_{0}$ , $r_{0}=0.5$.  It is worth noting that
the exact solution is $C_{0}$,  but $\partial u/\partial r$ is discontinuous
on the interface. 

The RBF distance network  is shown
in \Cref{fig:distance koch}b and c,  which has more collocation points than the crack
problem due to the complex boundary conditions.  We can find the
prediction of the RBF distance network matches the accurate distance. 

\begin{figure}
\begin{centering}
\includegraphics[scale=0.5]{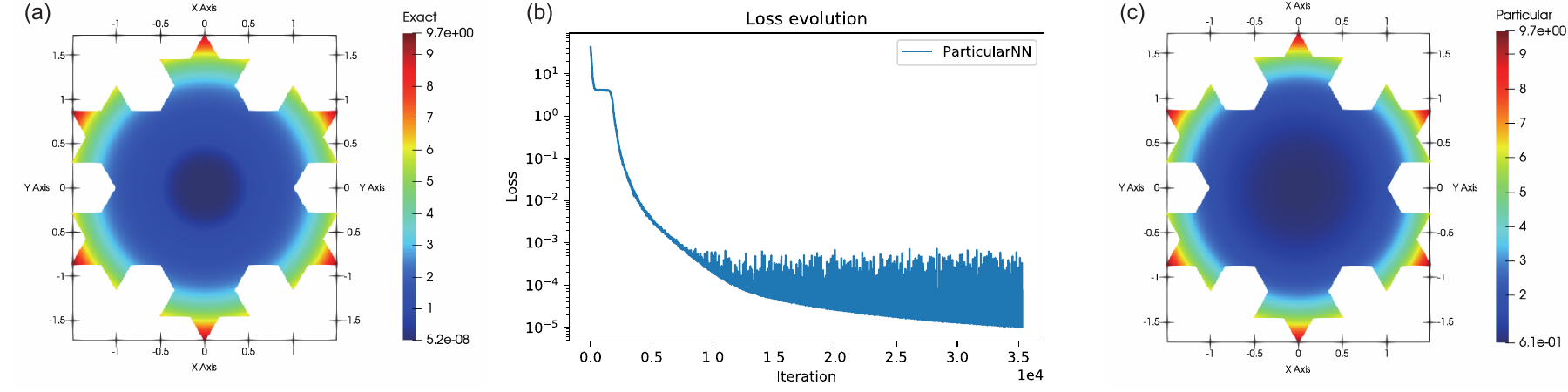}
\par\end{centering}
\caption{Koch particular neural network:  (a) The contour of the exact solution of
Koch. (b) The evolution of the MSE loss for the training particular
network. (c)  The contour of the  prediction of the particular network.\label{fig:particular-neural-network_koch}}
\end{figure}

The convergence condition of the particular network training is that
the MSE is less than 1e-6.  The hidden layer of the particular network
has 4 layers,  each layer has 20 neurons,  the activation function
is tanh in all hidden layers, and an identical function is used in the
output layer,  the optimizer is adam, the learning rate is 5e-4. \Cref{fig:particular-neural-network_koch}
shows the particular pattern is close to to the pattern of the exact solution, because the essential boundaries enclose the domain.  The loss
function value of the particular neural network decreases as the number
of iterations increases,  which shows the particular network can fit
the essential boundary accurately.  The oscillation of the loss is the reason
that the learning rate is not adequate small,  we can decrease the
learning rate as iterations increase to eliminate the oscillation phenomenon. 

Given that the derivative on the interface is discontinuous,
 we compare CPINN, DEM,  and CENN. We divide the domain into 2 subdomains, i.e., $r=0.5$ as the partition , considering the derivative discontinuity on the interface.  The neural network architecture is the same in the both subdomain,
 which can be also different.  The hidden layer of the particular
network has 4 layers,  each layer has 20 neurons.  The activation
function is tanh in all hidden layers, and identical function is used
in the output layer.  Adam is used as the optimizer, and the learning rate is 1e-3.
 The training points are randomly distributed in all methods as shown
in \Cref{fig:distance koch}a.  Training points are redistributed every 100 epochs for all three
methods.  The penalties $\lambda$ of every loss function in CPINN  are :
   $\lambda_{1}=60$ when $r<r_{0}$,  $\lambda_{2}=1$ when $r>r_{0}$,
 $\lambda_{inter}$=1 when $r=r_{0}$,  $\lambda_{bound}=30$ on the
essential boundary (We adjusted a lot of penalties $\lambda$, this is the best set).  It is worth noting that the loss of CPINN on the interface
includes 
\begin{equation}
	\begin{cases}
		u^{+}=u^{-}\\
		a_{1}(\boldsymbol{n}\cdot\nabla u^{+})=a_{2}(\boldsymbol{n}\cdot\nabla u^{-})
	\end{cases}
\end{equation}
where $u^{+}$ and $u^{-}$ represent the internal and external region
neural network respectively.  For the sake of simplicity,  the two interface
loss terms are controlled by the same penalty $\lambda_{inter}$.  We
adopt sequential training for two neural networks \citet{thermochemical_PINN_optimization_way}.
 In particular,  we first optimize $u^{+}$ parameters with $u^{-}$
parameters fixed for 20 epochs,  then vice versa in CPINN. DEM and
CENN are trained by the admissible function.  CENN is trained by
 sequential training to pay attention to the internal energy (internal),
 given that the difference between the energy of internal and external
is too large.  we first optimize $u^{+}$ parameters with $u^{-}$parameters
fixed for 2500 epoch in CENN,  then vice versa.  We use the admissible
function in the external region,  while not in the internal region 
 because the whole boundary condition is in the external region boundary.  We
stitch the subdomains through the interface loss. 

\begin{figure}
\begin{centering}
\includegraphics[scale=0.65]{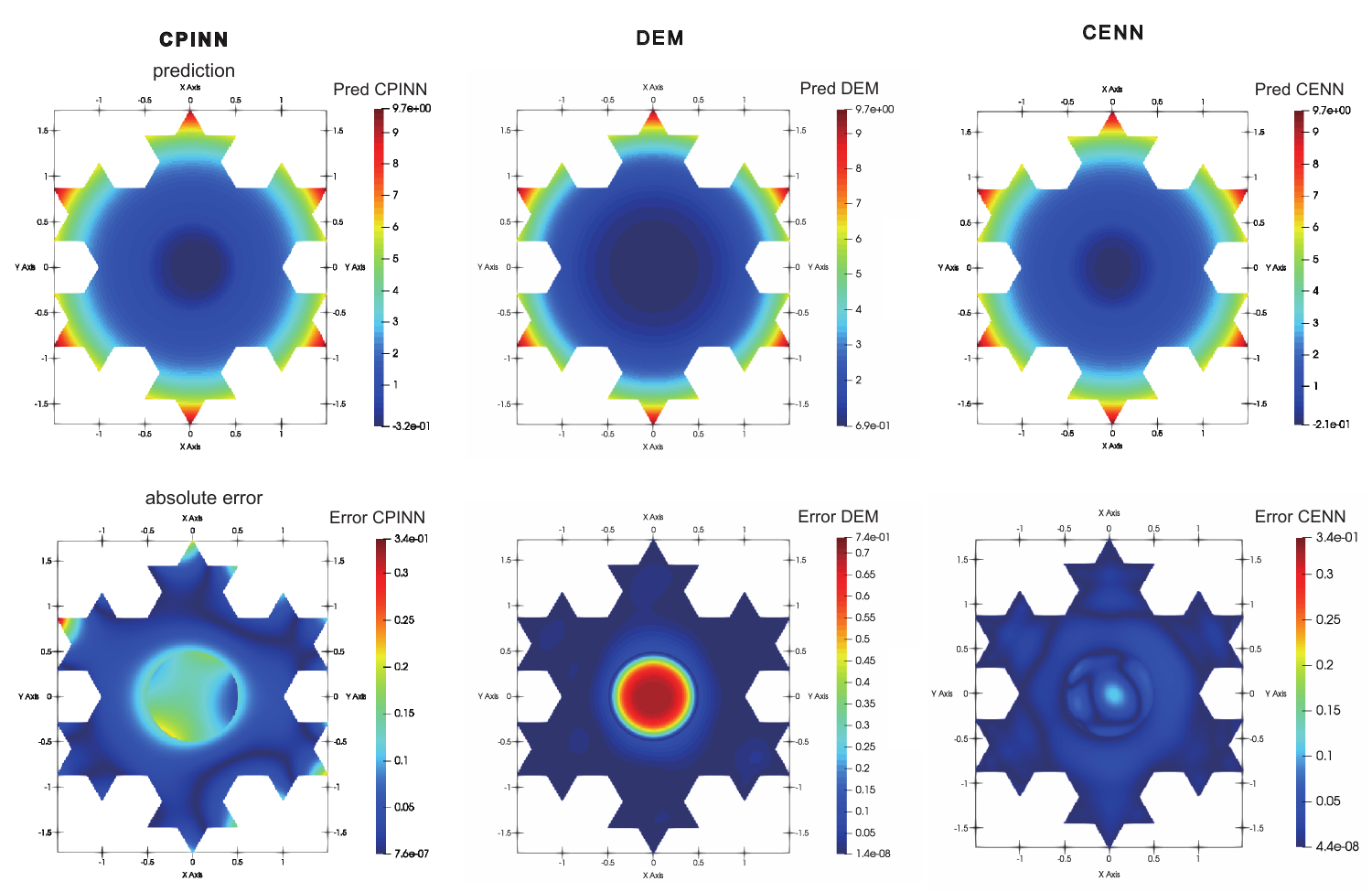}
\par\end{centering}
\caption{Prediction and absolute error by CPINN, DEM, and CENN in Koch: The
displacement field prediction result (a,b,c), absolute error (d,e,f). \label{fig:Koch comparison}}
\end{figure}

\Cref{fig:Koch comparison} shows the prediction and absolute error
of these three methods.  The maximum absolute error of CPINN is on
the boundary because  the second-order derivation
of the solution is big on the boundary in CPINN, resulting in the
prediction sensitivity on the boundary.  The error of DEM is mainly
in the middle of the region($r<r_{0}$),  because the energy density
$\frac{1}{2}\alpha_{1}(\nabla u)\cdot(\nabla u)=8r^{6}/a_{1}$ is
much smaller than the external energy density. The big difference in internal and external energies results in paying
too much attention to the external region. The maximum absolute error
in CENN, CPINN, CENN is 0.34, 0.74, and 0.11 respectively. The error
in the CENN is mainly concentrated in the area near the center point,
 which is obviously smaller than DEM in terms of error magnitude and
range. In CENN, because the energy density $\frac{1}{2}\alpha_{1}(\nabla u)\cdot(\nabla u)=8r^{6}/a_{1}$
is close to zero when $r$ is close to the center point, the loss
of the internal will pay attention to the region where the loss can
be decreased faster,  resulting in the error in the center point. 

\begin{figure}
\begin{centering}
\includegraphics[scale=0.52]{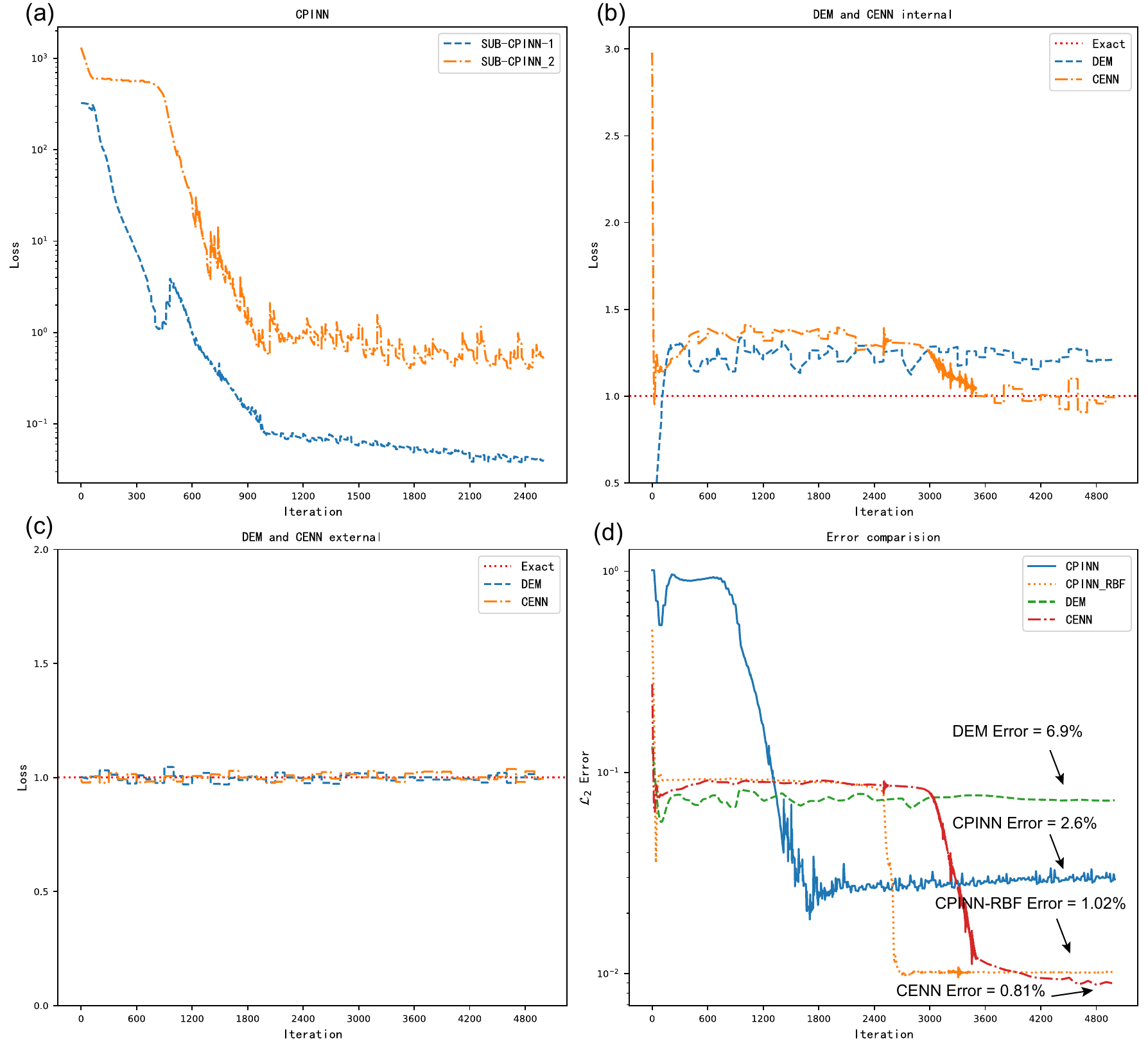}
\par\end{centering}
\caption{Comparison of CPINN, DEM, and CENN in Koch: (a) The loss function evolution
of CPINN. SUB-CPINN-1 and 2 are CPINN subdomains NN in the internal and external regions respectively. (b) Comparison of DEM and the CENN loss
function in the internal region with exact normalized internal functional
integration. (c) Comparison of DEM and the CENN loss
function in the external region with exact normalized functional
integration. (d) Comparison   of CPINN, CPINN-RBF, DEM, and CENN $L_2$ error.\label{fig:Koch error and loss}}
\end{figure}

To further analyze the performance of these three methods,  we compare
the loss and relative error $L_{2}$. Considering that the loss of
DEM and CENN are both energy variational,  so the optimal values both
are the same,  so we compare the loss functions of these two methods together
as shown in \Cref{fig:Koch error and loss}.  Due to the complexity
of the boundary conditions,  we calculate the exact numerical integral
for this case by the Monte Carlo algorithm. The exact values of the internal and the external are 1.10 and 436.45 respectively. To better visualize it,  we normalize the energy
loss.  \Cref{fig:Koch error and loss}a shows the trend
of the different part loss evolution in CPINN are similar,  which
proves our choice of the penalty is relatively good.  \Cref{fig:Koch error and loss}b shows DEM can not converge to the exact internal energy,
 which is the reason why the absolute error in the internal region
is obviously large than the external region.  However,  CENN can converge
to the exact internal energy well due to the subdomains.  \Cref{fig:Koch error and loss}c shows CENN and DEM both converge to the external energy 
 because the external energy dominates the majority of the total energy.
 The fluctuation in the external and internal energy can be eliminated
by decreasing the learning rate.  \Cref{fig:Koch error and loss}d shows the relative error $L_{2}$ is 2.6\%, 1.02\%,  6.9\%, and 0.81\% in CPINN, CPINN-RBF, DEM, and CENN respectively. We find that CPINN-RBF improves the performance against CPINN except in the first numerical experiment (crack).
It may be that the penalty factor in the first numerical example is more suitable than admissible function.  Obviously ,  CENN can solve
the composite problem very well,  and there is only one hyperparameter
penalty about the interface in the loss function.  Finally,  we can find that
the trends of the loss function and relative error of the above three
methods are the same. 

\begin{figure}
\begin{centering}
\includegraphics[scale=0.5]{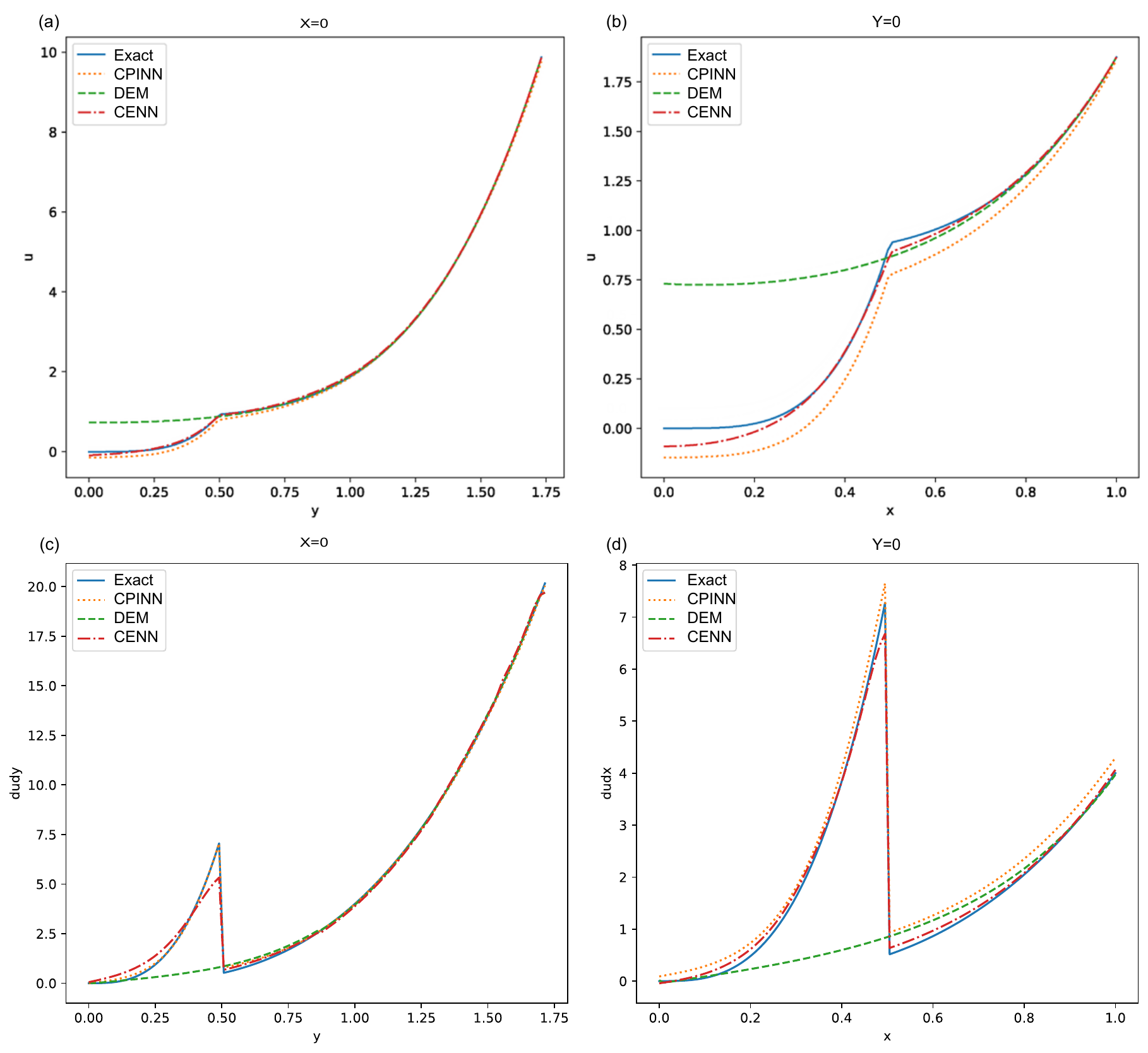}
\par\end{centering}
\caption{Different location comparison of CPINN, DEM, and CENN in the Koch: (a) Comparison
of the displacement solution on the vertical line x=0. (b)
Comparison of the displacement solution on the horizontal line y=0.
(c) Comparison of the displacement  derivative 
${\partial u}/{\partial r}$ on the vertical line x=0. (d)
Comparison of the displacement derivative ${\partial u}/{\partial r}$
on the horizontal line y=0. \label{fig:Different-location-koch}}
\end{figure}

To further analyze the interface prediction of these three methods,
 we compare the solution and the discontinuous derivative on the interface.
 \Cref{fig:Different-location-koch}a and b show  the comparison of
solution of CPINN, DEM and CENN at two cross lines,  x=0 and y=0.
CENN is more accurate than other methods in predicting solutions,
 and the error increases slightly near the center area,  the reason
may be relatively little energy density in the center point.  The
center error may be alleviated by adopting more points in the internal
region,  which we will further study in the future.  It is worth noting
that the derivative solution by CPINN is more accurate than CENN and
energy method,  as shown in \Cref{fig:Different-location-koch}c,d.
 The reason may be the second-order derivative used in CPINN is not
smaller than the first-order derivative in CENN in the internal region,
so the derivative in the internal region can be trained well by strong
form.  It is worth noting the solution prediction of these three methods
is smaller than the exact solution,  which is suitable for some engineering
conservative problems. 

Given that the strong form of PINN can adopt batch strategy, it is very useful to improve the efficiency of network training. However, because CENN is a deep energy method, the precise numerical integration of energy functional is very important. So we use as many integration points as possible in CENN. In order to compare the computational efficiency of CPINN and CENN, we compare CENN with CPINN under different batch sizes. \Cref{tab:Compared-time} shows that the efficiency and accuracy of CENN are better than CPINN in the same batch size 10000. The error and time in \Cref{tab:Compared-time} are the results after the network converges. We can find that different batch sizes of CPINN have a significant impact on the accuracy of the solution, which is  mainly due to the high randomness of the training points and the strong non-convex property of the loss function.
It is worth noting that the efficiency of CENN in solving high-order tensors and high-order PDEs will be greater than that of CPINN (the problem of crack is second-order and scalar PDEs, so the efficiency advantage of CENN is not very obvious ).

\begin{table}
	\caption{Compare  the training time of CPINN and CENN. The network structure
		is exactly the same. Batch size refers to the random number of points
		required for each epoch of CPINN, and each epoch is randomly scattered.
		CENN takes as many integration points as possible because of the accuracy
		requirement of energy numerical integration. $L_{2}$ error and time are the results of convergence. \label{tab:Compared-time}}	
	\begin{centering}
		\begin{tabular}{llll}
			\toprule 
			& Batch size & $\mathcal{L}_{2}$ error & Time (s)\tabularnewline
			\midrule
			& 100 & 7.1\% & 62.45\tabularnewline
			& 500 & 6.9\% & 66.15\tabularnewline
			\multirow{2}{*}{CPINN} & 1000 & 2.5\% & 71.25\tabularnewline
			& 3000 & 3.1\% & 176.25\tabularnewline
			& 5000 & 3.4\% & 148.2\tabularnewline
			& 10000 & 2.6\% & 247.8\tabularnewline
			\midrule 
			CENN (Whole Batch) & 10000 & 0.81\% & 176.4\tabularnewline
			\bottomrule
		\end{tabular}
		\par\end{centering}

\end{table}

\subsection{Composite material hyperelasticity\label{sec:hpyerelasticity material}}

In this section,  we introduce hyperelasticity,  which is a well-known 
problem including non-linear operator and vector-valued variables
in solid mechanics.  In \citet{PINN_hyperelasticity},  the homogeneous
hyperelasticity problem is solved.  We consider a body made of a nonhomogenous material to show the proposed method CENN can tackle the non-linear and heterogeneous problem,
 i.e., composite hyperelastic material.  The non-linear characteristic
is from the large deformation and the hyperelastic constitutive law.
 The governing equation of the problem is 
\begin{equation}
	\begin{cases}
		\nabla_{X}\cdot \boldsymbol{P}+\boldsymbol{f}=0 & \boldsymbol{x}\in\Omega\\
		\boldsymbol{u}=\bar{\boldsymbol{u}} & \boldsymbol{x}\in\partial\Omega^{eb}\\
		\boldsymbol{N}\cdot\boldsymbol{P}=\bar{\boldsymbol{t}} & \boldsymbol{x}\in\partial\Omega^{nb}
	\end{cases}\label{eq:elastic strong form}
\end{equation}
where $\nabla_{X}$ is the gradient operator with respect to $\boldsymbol{X}$,
$\boldsymbol{X}$ is material coordinate \citet{the_foundation_of_solid_mechanics_feng}.
 $\nabla_{X}\cdot \boldsymbol{P}$ denotes the divergence operator,  and we can
use the tensor index to represent clearly $\nabla_{X}\cdot \boldsymbol{P}=P_{ij,i}$.
 $\boldsymbol{P}$ is the Lagranges's stress,  where the first component $i$ and
the second $j$ of $P_{ij}$ correspond the material coordinates and
spatial coordinate respectively.  $\boldsymbol{f}$ is the body force,  and the
first equation is the equilibrium equation in the domain $\Omega$.
 Note that $\boldsymbol{P}$ is the function of the material coordinate $\boldsymbol{X}$.  $\bar{\boldsymbol{u}}$
is the given displacement value in the essential boundary $\partial\Omega^{eb}$.
 $N$ is the normal direction of the Newmann boundary $\partial\Omega^{nb}$,
 and $\bar{\boldsymbol{t}}$ is the prescribed traction on the Newmann boundary. 

If the material is hyperelasticity,  then $\boldsymbol{P}$ can be written to be the derivative
of the strain energy $\Psi$ with respect to the deformation gradient
$\boldsymbol{F}$
\begin{equation}
	\boldsymbol{P}=(\frac{\partial\Psi}{\partial\boldsymbol{F}})^{T}
\end{equation}
\begin{equation}
	\boldsymbol{F}=\frac{\partial\boldsymbol{x}}{\partial\boldsymbol{X}}
\end{equation}
where the $\boldsymbol{x}=\boldsymbol{X}+\boldsymbol{u}(\boldsymbol{X})$ is the spatial coordinate,  which is the function
of the material coordinate $\boldsymbol{X}$,  given that it is a static problem
independent of time.  $\boldsymbol{u}$ is the interesting field,  which in this
specific problem is the vector displacement field.  We consider the
common used constitutive law Neo-Hookean of the hyperelasticity \citet{belytschko2013nonlinear}
\begin{equation}
	\varPsi=\frac{1}{2}\lambda(lnJ)^{2}-\mu ln(J)+\frac{1}{2}\mu(trace(\boldsymbol{C})-3)
\end{equation}
where $J$ is the determinant of the deformation gradient $\boldsymbol{F}$ and
$\boldsymbol{C}$ is the Green tensor, i.e., $\boldsymbol{C}=\boldsymbol{F}^{T}\cdot \boldsymbol{F}$.  The first and second
terms of RHS respectively are uncompressible condition and stress-free
in the initial condition.  $\lambda$ and $\mu$ is Lame parameter
\begin{equation}
\begin{cases}
\lambda=\frac{vE}{(1+v)(1-2v)}\\
\mu=\frac{E}{2(1+v)}
\end{cases}
\end{equation}
where $E$ and $v$ represent the elastic modulus and Poisson ratio
respectively.  The core of the problem is to obtain the displacement
field $\boldsymbol{u}$.  There are two ways to do it.  The first way is to solve
the strong form \Cref{eq:elastic strong form} and the second way is
to use the energy theory by optimizing the potential energy

\begin{equation}
	\mathcal{L}=\int_{\Omega}(\Psi-\boldsymbol{f}\cdot\boldsymbol{u})dV-\int_{\partial\Omega^{nb}}\bar{\boldsymbol{t}}\cdot\boldsymbol{u}dA
\end{equation}

 Noting that the trial function satisfies the essential boundary in
advance.  In \ref{sec:Appendix C hyper strong}, we can find the implementation
of strong form is more complicated than the energy form, and the highest
derivative of the strong form  is higher than the energy form (The best hyperparameter combination  of different PDEs and boundary conditions in CPINN is not known), which results in more computation
cost and lower accuracy.  For the sake of simplicity,  we compare
DEM and CENN.  We use the FEM as the reference solution.

\begin{figure}
\begin{centering}
\includegraphics[scale=0.8]{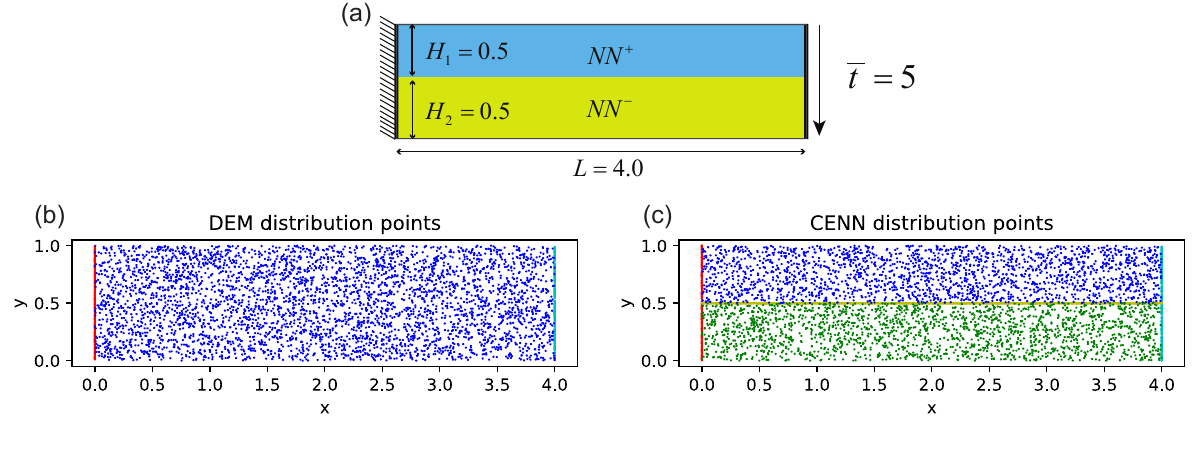}
\par\end{centering}
\caption{The schematic diagram of the composite material beam, the blue and
yellow stack denotes the different hyperelastic material (up: $E_{0}=1000,\mu_{0}=0.3$;
down: $E_{1}=10000,\mu_{0}=0.3$), we use different neural networks
corresponding to the different material. The length of the beam is
4 and the height of the beam is 1, the downward traction acting the
right end of the beam is 5N and the left end of the beam is an essential
fixed boundary (a). Illustration of distribution points strategy
of the DEM (b) and CENN (c), In DEM, the blue points
are the domain training points of the neural network. In CENN, The
blue points are the upper material training points of the neural network
$NN^{+}$, the green points are the lower material training points
of the neural network $NN^{-}$, the red points are the points of
the essential boundary, and the yellow points are the interface points
of the different materials. In both two energy ways, the total number
of internal points is both 4096, the number of essential boundary points
is both 256, and the number of interface points is both 1000. All points are
randomly distributed (b,c)\label{fig:The-schematic hyper}}
\end{figure}

\begin{figure}
\begin{centering}
\includegraphics[scale=0.8]{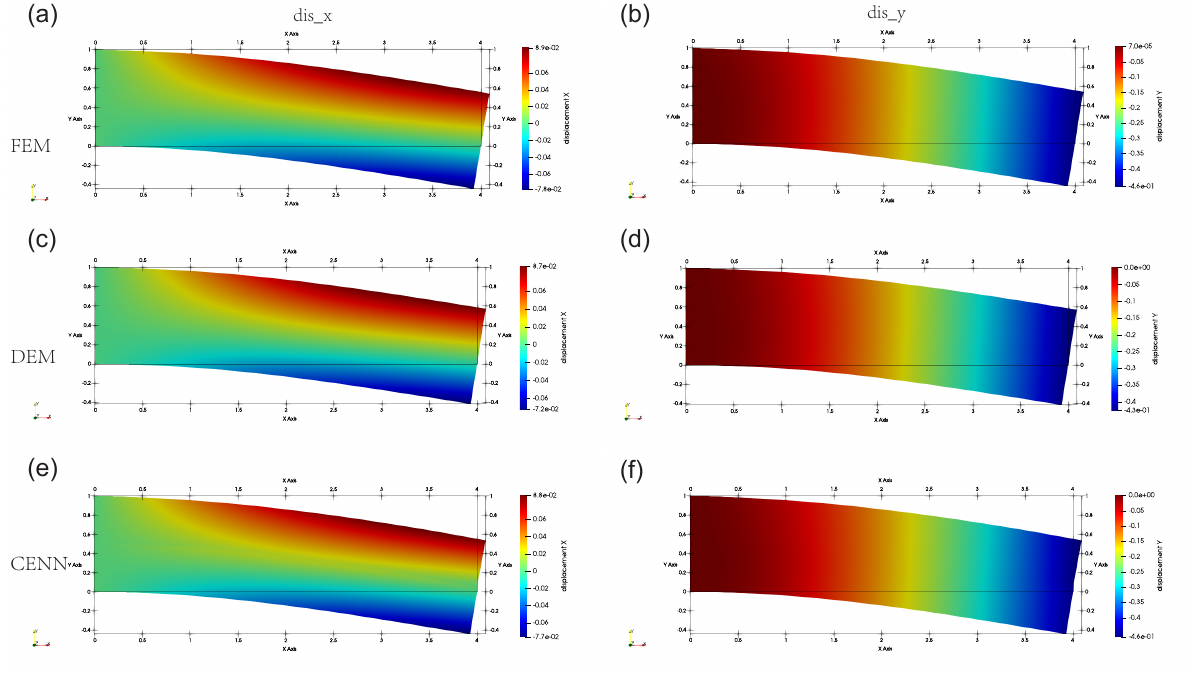}
\par\end{centering}
\caption{Predicted solution by FEM, DEM, and CENN in beam: The displacement
X field prediction result (a,c,e), the displacement Y field prediction
result (b,d,f), FEM predition (a,b), DEM predition (c,d),
CENN predition (e,f). \label{fig:predicted-soluton-beam}}
\end{figure}

\begin{figure}
\begin{centering}
\includegraphics[scale=0.5]{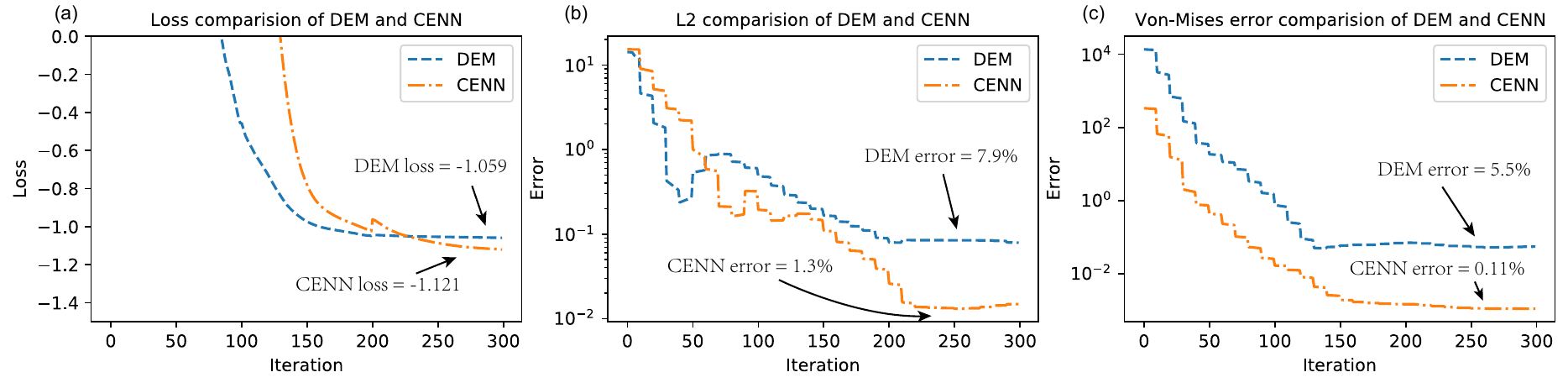}
\par\end{centering}
\caption{Comparison of DEM and CENN in beam: (a) The loss function evolution.
(b) The relative error $L_2$ about displacement magnitude. (c) The relative
error about Von-Mises stress. \label{fig:comparison-error_beam}}
\end{figure}

We consider a two-dimensional bending beam made of the composite material,
as shown in \Cref{fig:The-schematic hyper}. The material parameter
is, upper : $E_{0}=1000,\mu_{0}=0.3$, lower : $E_{1}=10000,\mu_{0}=0.3$.
FEM divides the region
into 400 (length){*}100 (height) quadratic elements for enough accurate
reference solution. \Cref{fig:The-schematic hyper} shows the points
distribution ways of the DEM and CENN.  Training points are redistributed
every 100 epochs in all methods.  Because the given value on the essential
boundary is zero,  the particular neural network is zero, i.e.,
$\boldsymbol{u}_{p}=0$.  In addition, given that the essential boundary geometry
is simple, we can obtain the analytic solution of the distance function, i.e., $RBF(x,y)=x$ to ensure that the boundary conditions are exactly
satisfied (RBF method can also be used here, but it is not necessary due to the simplicity of the boundary).  The structure of the generalized
network is 4 layers,  each layer has 20 neurons,  and the learning
rate is 1e-3.  We use LBFGS \citet{goodfellow2016deep} as the optimizer.
 Given that the strain energy does not vary much from region to region,
we optimize both neural networks in CENN not like the Koch example
before.  The neural network structure of DEM is same as CENN
for comparison.  \Cref{fig:predicted-soluton-beam} shows the prediction
of the FEM,  DEM, and the CENN.  We can find that the minimum displacement
x and y of the CENN is more close to the reference solution, i.e.
FEM.  The pattern of DEM and CENN both coincides with FEM.  This shows
the prediction of DEM and CENN about original field $\boldsymbol{u}$ is both good
overall.  To further quantify the error,  we consider the relative
error about displacement and Von-Mises stress related to the derivative
of the displacement
\begin{equation}
	Mises=\sqrt{\frac{3}{2}\boldsymbol{S}^{dev}:\boldsymbol{S}^{dev}}
\end{equation}
where $\boldsymbol{S}$ is the Kirchhoff's stress 
\begin{equation}
	\boldsymbol{S}=\boldsymbol{P}\cdot\boldsymbol{F}^{-T}.
\end{equation}
and $\boldsymbol{S}^{dev}$ is the deviation tensor, i.e., 
\begin{equation}
	\boldsymbol{S}^{dev}=\boldsymbol{S}-\frac{1}{3}trace(\boldsymbol{S})\boldsymbol{I}
\end{equation}
Because of the different material parameters of the beam,  the derivative
of the displacement y with respect to the y-direction is discontinuous
on the interface.  This is also the case with the Von-Mises  on
the interface,  so we consider the error $abs(\Vert Von^{pred}\Vert_{2}-\Vert Von^{fem}\Vert_{2})/\Vert Von^{fem}\Vert_{2}$
to consider the derivative error.  \Cref{fig:comparison-error_beam}a shows the loss evolution of the CENN is much lower than the
CENN.  It is worth noting the optimal loss of the energy method is
not zero due to  the minimal potential theory. In some ways,  the
loss can reflect the accuracy of the method.  \Cref{fig:comparison-error_beam}b shows the relative error of CENN about displacement magnitude
is 1.3\%,  which is lower than DEM. Because $\partial u/\partial y$
is not continuity, it is natural to divide the region according to
the different materials to ease the restriction of derivative continuity. \Cref{fig:comparison-error_beam}c shows the relative error of CENN in Von-Mises stress is
0.11\%,  which is very much smaller than DEM.  Because CENN can fit better
the derivative on the interface due to the inherent ability to
simulate the discontinuity,  the prediction about Von-Mises in CENN
is more close to the reference solution. 

\begin{figure}
\begin{centering}
\includegraphics[scale=0.5]{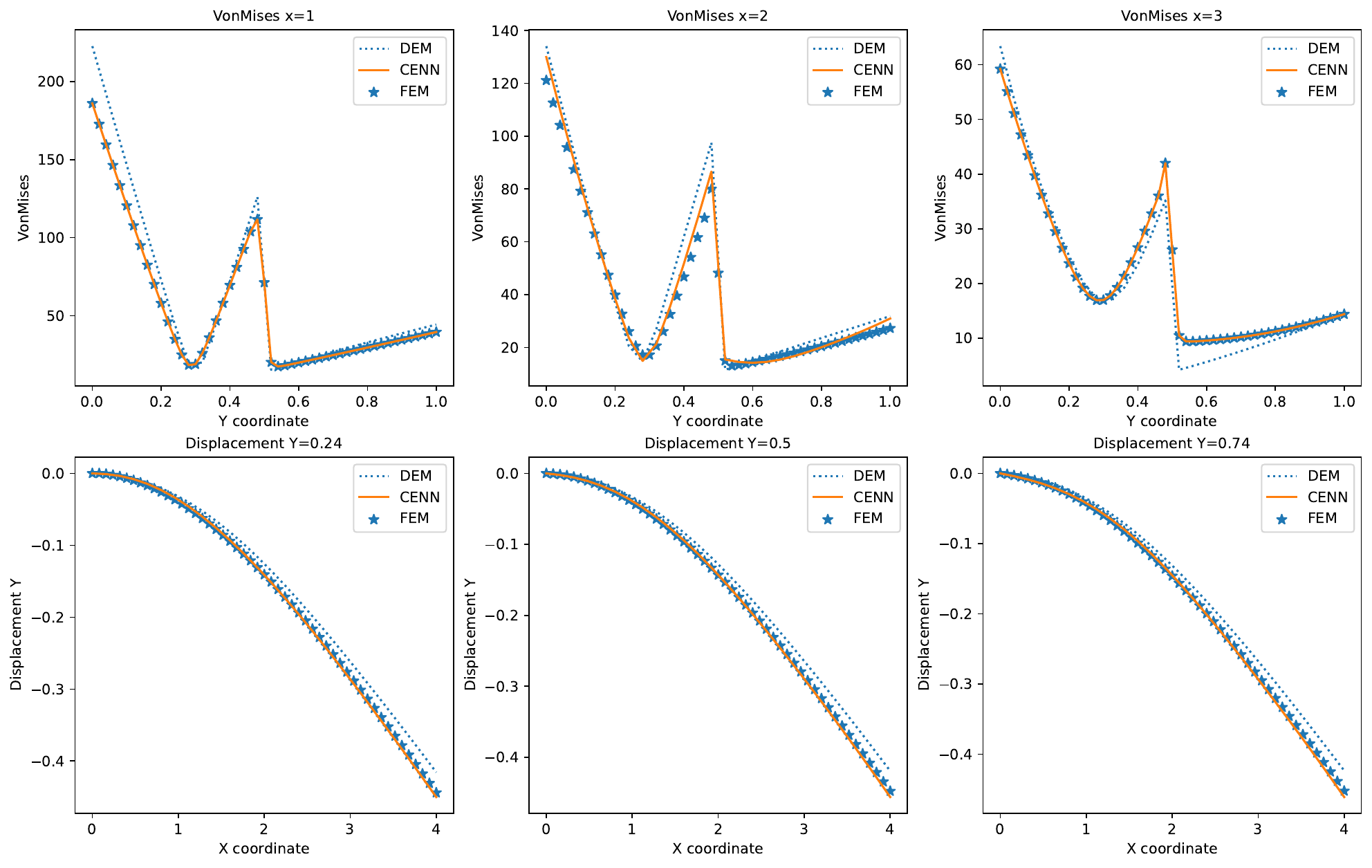}
\par\end{centering}
\caption{Different location comparison of DEM and CENN in beam: Comparison
of the Von-Mises solution on the different vertical line, x=1, x=2,
x=3 respectively (a,b,c), comparison of the displacement y solution
on the different horizontal line, y=0.25, y=0.5, x=0.75 respectively
(d,e,f). \label{fig:Different-location-comparion_beam}}
\end{figure}

To further analyze the interface prediction of DEM and CENN,  we compare
the discontinuous derivative on the different vertical lines cross
the interface,  x=1,  x=2, and x=3,  and we also compare the main deformation,
 i.e. $u_{y}$,  on the different horizontal line,  y=0.25,
 y=0.5 and y=0.75. \Cref{fig:Different-location-comparion_beam}a,b,c
show the comparison of Von-Mises of DEM and CENN on the different
vertical lines cross the interface,  x=1,  x=2, and x=3.  We can find
CENN is more accurate than DEM in terms of discontinuous derivative.
 \Cref{fig:Different-location-comparion_beam}d,e,f show the comparison
of displacement y of DEM and CENN on the different horizontal lines,
 y=0.25, y=0.5, and y=0.75.  We can find both methods are good overall
in terms of displacement y,  but CENN is better.  

\section{Discussion\label{sec:Discussion}}

\subsection{Discussion of admissible function errors}

The above-mentioned data-driven and CENN error at the center point (crack problem) is
larger than other points.  Here we analyze the reason.  We consider
the derivative of the ring direction $\theta$ at the center point.
\begin{equation}
	\frac{\partial u(r,\theta;\boldsymbol{\theta}_{p},\boldsymbol{\theta}_{g})}{r\partial\theta}=\frac{\partial u_{p}(r,\theta;\boldsymbol{\theta}_{p})}{r\partial\theta}+\frac{\partial RBF(r,\theta)}{r\partial\theta}\cdot u_{g}(r,\theta;\boldsymbol{\theta}_{g})+RBF(r,\theta)\cdot \frac{\partial u_{g}(r,\theta;\boldsymbol{\theta}_{g})}{r\partial\theta}\label{eq:Distance network loop derivative}
\end{equation}

\begin{figure}
\begin{centering}
\includegraphics{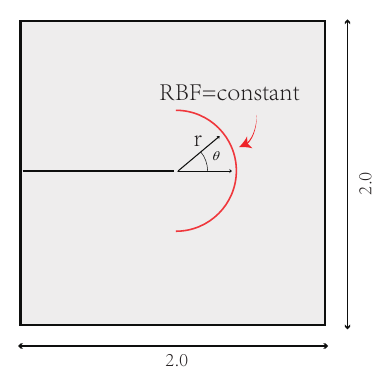}
\par\end{centering}
\caption{RBF distance network constant diagram, r and $\theta$ are standard
polar coordinates, the prediction of RBF distance network is constant
in the right half when $r=constant$, and the red line is the position
where RBF prediction is constant\label{fig:Schematic diagram of :RBF distance function constant}}
\end{figure}

Here we assume that the RBF distance network is correct.  It can be
found in \Cref{fig:RBF}b ,  the result of the RBF distance
network is consistent with the analytical solution.  Therefore,  this
assumption is reasonable.  Obviously,  the RBF distance network in
the \Cref{eq:Distance network loop derivative} is equal to zero at
the center point.  At the same time, the RBF distance network keeps
the same distance from the center point on the right half,  as shown
in \Cref{fig:Schematic diagram of :RBF distance function constant},
i.e., $\partial RBF(r,\theta)/r\partial\theta=0$. So when approaching
the center point, $\partial RBF(r,\theta)/r\partial\theta$ and $RBF(r,\theta)$
are both zero,
which will cause $\partial u(r,\theta)/r\partial\theta=\partial u_{p}(r,\theta)/r\partial\theta$.
However we fixed the parameters of the particular network during training,
so the energy density $J_{\rho}=\frac{1}{2}[(\partial u/\partial r)^2+(\partial u/r\partial\theta)^{2}]$
at the center point can not be trained completely, i.e.,  $\partial u/r\partial\theta$
does not change,  which leads to a large  error at the center.  \Cref{fig:Comparison_predition_error_crack}e and f show that the error is mainly concentrated
in the right half of the center point.  This shortcoming can be overcomed
by changing the derivative term of the RBF distance network, i.e.,
change the form of the distance network to make $\partial RBF(r,\theta)/r\partial\theta\neq0$
at the center point. It will make training of the generalized network
useful to perform learning further in the center. In another way, we
can retrain the parameters of the particular network, but we need
to add an additional penalty function to limit the change of the particular
network parameters so that the energy density at the center point
can be changed. This problem will be further studied in the future. 

We consider the singular strain at x>0,  y=0.  \Cref{fig:Comparison of different locations of cracks}d shows that  the error will increase when approaching x=0.
 This is actually a special case of the above-mentioned hoop $\theta$
derivative situation at the interface (x>0, y=0), $\partial u(r,\theta)/r\partial\theta=\partial u(x\text{,}y)/\partial y$.
 It is also caused by the RBF distance network of the addmissible
function,  we consider the derivative $u$ w.r.t y-direction at the
center point 
\begin{equation}
\frac{\partial u(x\text{,}y;\theta_{p},\theta_{g})}{\partial y}=\frac{\partial u_{p}(x,y;\theta_{p})}{\partial y}+\frac{\partial RBF(x,y)}{\partial y}\cdot u_{g}(x,y;\theta_{g})+RBF(x,y)\cdot \frac{\partial u_{g}(x,y;\theta_{g})}{\partial y}.\label{eq:Possible displacement field derivatives}
\end{equation}
Obviously,  the RBF distance network in \Cref{eq:Possible displacement field derivatives}
is equal to zero at the center point.  In addition, we analyze $x=\delta$,
 y=0 at $\partial RBF(x,y)/\partial y$,  Taylor series expand to
second-order term
\begin{equation}
\begin{split}\frac{\partial R(y)}{\partial y} & =lim_{d\rightarrow0}\frac{\sqrt{\delta^{2}+d^{2}}-\delta}{d}\\
 & =lim_{d\rightarrow0}\frac{(\sqrt{\delta^{2}}+\frac{\partial(\sqrt{\delta^{2}+d^{2}})}{\partial y}|_{d=0}d+\frac{1}{2}\frac{\partial^{2}(\sqrt{\delta^{2}+d^{2}})}{\partial^{2}y}|_{\triangle=0}d^{2}+o(d^{2}))-\delta}{d}\\
 & =lim_{d\rightarrow0}\frac{(\sqrt{\delta^{2}}+\frac{1}{2}\frac{1}{\delta}d^{2}+o(d^{2}))-\delta}{d}\\
 & =0
\end{split}
\end{equation}
Where $\delta$ is the distance from the center point (crack tip),
 $\delta>0$,  the coordinate point is $(x=\delta,y=0)$,  where $d$
is the distance from the coordinate point to the x-axis,  it is not
difficult to find that the above derivative is always equal to 0
at $x>0,y=0$. So when approaching to the center point,  $\partial R(x,y)/\partial y=R(x,y)=0$,
 which will lead to $\partial u(x\text{,}y)/\partial y=\partial u_{p}(x,y)/\partial y$. However,  we fixed the parameters of the particular network during training,
 so the closer is to the center point,  the more sensitive generalized
network is. It will increase  the error of $\partial u(x\text{,}y)/\partial y$ when approaching a center point. 

\subsection{The influence of point distribution on the solution}

Different point allocation methods may affect the accuracy.  At the
same time,  due to the constant of the distance network near the essential
boundary,  the generalized neural network near the essential boundary
partially is very sensitive, resulting in increasing error. The
sensitivity can be reduced by arranging more points around the center
point. 

The energy form depends more on the integration scheme because the loss function is the value of the energy functional. The different ways of attribution points have a relatively large impact on the energy form. However, the strong form requires the error at the sampling point to be 0, so the requirement for the integration scheme is not high. It is worth noting that the overall loss function tends to be 0 in strong form, so the gradient descent has an impact to all the sampling points. In a way, the attribution ways have also an impact on the strong form. The impact of the attribution ways on the energy form is shown explicitly, while the impact of the attribution ways on the strong form is implicit.
\subsection{Efficiency and accuracy}

Since the addmissible function requires additional training,  the
advanced training of the RBF distance network and the particular network
will increase the additional computational cost.  However,  the traditional
strong form does not need it,  but the traditional strong form involves
a higher-order derivative than the energy form. Thus the computational
cost of the strong form is larger than the energy form after the admissible
function training.  This is a question of balance.  In addition,
 higher-order derivatives will not only further increase the amount
of the computational cost,  but also affect the accuracy theoretically.
It is worth noting that not all PDEs have a corresponding energy form,
 so the generality of the energy form is not as good as the strong
form.  At the same time,  the energy form requires more strict
collocation methods due to the requirement for precise integration functionals,
 while the strong form does not have this restriction,  and some batch
methods can be used to improve convergence speed \citet{the_comparision_of_strong_and_energy_form}. In addition,  the choice of hyperparameters in traditional strong
form is also a problem.  Although there is NTK theory to automatically
select hyperparameters \citet{NTK_PINN,NTK_to_get_hyperparameter_of_PINN},
 accurate and quantitative optimal hyperparameters is still an
important problem.  Of course,  the strong form can also use the construction
of the addmissible function to reduce the number of hyperparameters.
 If multiple neural network subdomains are involved, the hyperparameters
about the interface will further increase greatly,  and the interface
loss term in CPINN will add the derivative term compared to the CENN.
 It will reduce the accuracy and efficiency.  It is worth noting that
there is no NTK theory of subdomains to automatically determine the
hyperparameters,  which will be further studied in the future. 

\subsection{Use adaptive activation function to accelerate the convergence and improve the accuracy}
There are some adaptive activation methods to accelerate the convergence
and improve the accuracy  \citet{KNN_adaptive, LAAF,
Adaptive_activation_functions}. The main motivation of the adaptive
activation function is to increase the slope of the activation function
 \citet{LAAF}. The advantage of these methods is easy
and useful, and we only add a few trainable parameters to the original
neural network (the extra parameters are far less than the initial
trainable parameters). We include LAAF\citet{LAAF}
and Rowdy\citet{KNN_adaptive} to our initial CENN.
LAAF adds extra trainable parameters in every layer of the neural
network. The modification of LAAF compared to the origin neural network
is
\begin{equation}
\sigma(na^{(l)}\boldsymbol{z}^{(l)})
\end{equation}
, where untrainable n\ensuremath{\ge} 1 is a pre-defined scaling factor
to control the convergence speed of L\_LAAF and the parameter trainable
$a$ acts as a slope of activation function. The loss term $\mathcal{L}_{a}$
about the $a$ need to be added to the loss function, 
\begin{equation}
	\begin{aligned}\mathcal{L}_{LAAF} & =\mathcal{L}_{initial}+\lambda_{a}\mathcal{L}_{a}\\
		\mathcal{L}_{a} & =\frac{D}{\sum_{k=1}^{D}exp(a^{(k)})}
	\end{aligned}
\end{equation}
, where $D$ is the number of the hidden layers (not including input
and output layers). The meaning of $\mathcal{L}_{a}$ is the reciprocal
of the mean $exp(a^{(k)})$ in every hidden layer. $\lambda_{a}$
is the weight of $\mathcal{L}_{a}$, we adopt $\lambda_{a}=1$. The
detail of the L\_LAAF is in \citet{LAAF}.

Rowdy activation function is the special form of the Deep Kronecker
neural networks (A general framework for neural networks with adaptive
activation functions), which is demonstrated good performance in some
numerical experiments \citet{KNN_adaptive, LAAF, Adaptive_activation_functions}. The
modification of Rowdy compared to the origin neural network is
\begin{equation}
	\sum_{k=1}^{K}c_{k}^{(l)}\sigma_{k}(a_{k}^{(l)}\boldsymbol{z}^{(l)})
\end{equation}
, where $c_{k}^{(l)}$ and $a_{k}^{(l)}$ are the trainable parameters
in every hidden layer (the number of these is $2KD$). K is the different
activation function. The Rowdy activation function which is a special
different activation function is 
\begin{equation}
\sigma_{k}(x)=nsin((k-1)nx)
\end{equation}
where n \ensuremath{\ge} 1 is the fixed positive number acting as
a scaling factor. The detail of the Rowdy activation function is in
\citet{KNN_adaptive}.

\begin{figure}
	\begin{centering}
		\includegraphics[scale=0.5]{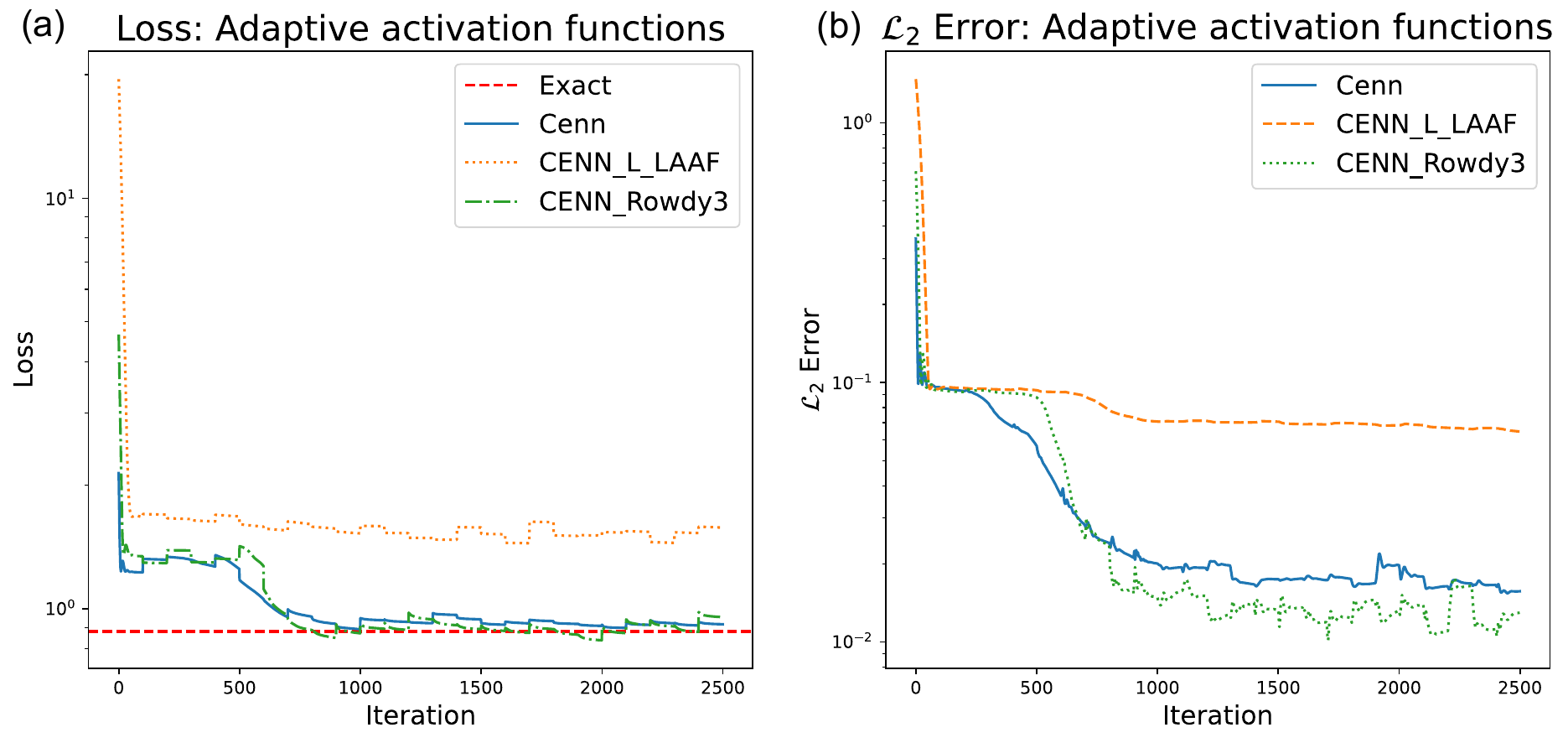}
		\par\end{centering}
	\caption{The performace of the adaptive activation function in crack example:
		Comparison of CENN, CENN with L\_LAAF, and CENN with Rowdy3. (a) Loss
		function (b) Total relative $\mathcal{L}_{2}$ error.\label{fig:Activation function}}
\end{figure}

We test LAAF and Rowdy3 in our first numerical example (crack), to
examine the performance of the adaptive activation function in CENN.
There are 3 different activation functions in Rowdy3, 
\begin{equation}
\begin{aligned}\sigma_{1}(x) & =tanh(x)\\
	\sigma_{2}(x) & =sin(x)\\
	\sigma_{3}(x) & =sin(2x)
\end{aligned}
\end{equation}
K=3, n=1 in Rowdy. n=10 in L\_LAAF. The architecture of the neural
network is the same as the \Cref{sec:crack problem}. The initialization of $a^{(l)}$
is 0.1 in all hidden layers in L\_LAAF. The initialization of $c_{k}^{(l)}$
is $\left[\begin{array}{cccc}
	1 & 0 & \cdots & 0\end{array}\right]_{k*1}^{(l)}$ in every hidden layer. The initialization of $a_{k}^{(l)}$ is $\left[\begin{array}{cccc}
	1 & 1 & \cdots & 1\end{array}\right]_{k*1}^{(l)}$ in every hidden layer. \Cref{fig:Activation function} shows CENN
with Rowdy3 can accelerate the convergence and improve the accuracy
compared to the initial CENN without adaptive activation function.
However, the CENN with L\_LAAF does not improve accuracy and convergence
speed. The easy numerical example shows that we can use the adaptive
activation function to accelerate convergence and improve accuracy
in CENN.

\subsection{Whether the generalized neural network is slicing}

Considering that the interface is a continuous derivative in the crack,
we do not actually need to slice the generalized network. But in some
$C_{0}$ problems (the original function is continuous,  the first
derivative is discontinuous), e.g. heterogeneous problem, the slicing
of the generalized network is meaningful,  which has been specifically
illustrated in \Cref{sec:non homogeneous boundary problem} and \Cref{sec:hpyerelasticity material}. 

We should consider the physical nature of the problem, then decide whether to slice the neural network or not. The different regions are connected through the interface loss function.

\section{Conclusion\label{sec:Conclusion}}

We developed the CENN formulation based on the deep energy method
to solve the heterogeneous PDEs. We showed
the accuracy of the proposed method is higher than the DEM in terms
of the heterogeneous problem,  especially the derivative on the interface.
 We show the efficiency and accuracy of the proposed method CENN compared
to the CPINN due to the lower-order derivative. The hyperparameter
of the proposed method is far less than CPINN.  We show CENN can solve
the singularity strain problem,  i.e. crack, without special treatment
such as finite crack element. We demonstrate why the PDE solution that
uses neural networks as approximation functions can be successful
from the perspective of trial functions as well as test functions.
We proposed a method to construct the admissible function to solve
the complex boundary problem,  which extends the application range
of the deep energy method.  Further,  we explore the  ability of the proposed method 
to solve the hyperelasticity problem with heterogeneous material,
 high-order multi-physics, and vector-valued variables in solid mechanics.
We show the ease of implementation of CENN compared to CPINN. We recommend
the CENN if the heterogeneous PDE has the corresponding variational
form because the accuracy and efficiency are better than CPINN and
DEM.

Each region can be set individually in CENN,  such as neural network
structure,  optimization method,  this flexibility is a
double-edged sword,  the advantage is the freedom of choice,  the
disadvantage is that there is no optimal solution.  A limitation of
the proposed method is that the derivative w.r.t some direction on
the essential boundary can not be learned well,  which we will further
study in the future. A limitation of this study is that the penalty
of the interface is constructed by a heuristic algorithm,  so the
accurate quantity of the hyperparameter should be further studied
in the future.  CENN is based on the principle of the minimum potential energy, so it only can solve the static problem. In the future, we will research the Hamilton principle to solve the dynamic problem in PINN energy form. An additional uncontrolled factor in almost scientific
computation based deep learning is how to choose initial parameters,
optimization way,  neural network architecture, and so on. It is worth emphasizing that the method has a natural advantage in dealing with heterogeneous problems.

\section*{Acknowledgement}
The study was supported by the Major Project of the National Natural Science Foundation of China (12090030).  The authors would like to thank Chenxing Li, Yanning Yu, and Vien Minh Nguyen-Thanh  for helpful discussions.

\appendix

\section{Back propogation of CENN\label{sec:Appendix-B.-Back propogation of CENN}}

We consider the energy form of CENN 
\begin{equation}
	\begin{split}\mathcal{L} & =\sum_{i=1}^{n_{1}}w_{i}w_{\varepsilon}(u_{1}(\boldsymbol{x}_{i};\theta_{1}))-\sum_{i=1}^{n_{1}}w_{i}f_{i}u_{1}(\boldsymbol{x}_{i};\theta_{1})\\
		& +\sum_{i=1}^{n_{2}}w_{i}w_{\varepsilon}(u_{2}(\boldsymbol{x}_{i};\theta_{2}))-\sum_{i=1}^{n_{2}}w_{i}f_{i}u_{2}(\boldsymbol{x}_{i};\theta_{2})\\
		& +\frac{penalty}{n_{inter}}\sum_{i=1}^{n_{inter}}(u_{1}(\boldsymbol{x}_{i};\theta_{1})-u_{2}(\boldsymbol{x}_{j};\theta_{2}))^{2},
	\end{split}
	\label{eq:CENN loss}
\end{equation}
where $w_{i}$ is the weight of the attribution points $x_{i}$,  especially
$w_{i}=V/n$ if uniform random Monte Carlo method is adopted. $w_{\varepsilon}$
is the strain energy density, $u_{1}$and $u_{2}$ are interesting
variable according to the different subdomains. $f_{i}$ is the body
force w.r.t. the point $x_{i}$. $\theta_{1}$and $\theta_{2}$ are
neural network parameters according to the different subdomains. The
last term on RHS is the continuity condition on the interface by the MSE
criterion.

Back propagate \Cref{eq:CENN loss} to find the derivative of the loss
function with respect to parameters of neural network 1
\begin{equation}
	\begin{split}\frac{\partial\mathcal{L}}{\partial\theta_{1}} & =\sum_{i=1}^{n_{1}}w_{i}\frac{\partial w_{\varepsilon}(u_{1}(\boldsymbol{x}_{i};\theta_{1}))}{\partial u_{1}(x_{i};\theta_{1})}\frac{\partial u_{1}(\boldsymbol{x}_{i};\theta_{1})}{\partial\theta_{1}}-\sum_{i=1}^{n_{1}}w_{i}f_{i}\frac{\partial u_{1}(\boldsymbol{x}_{i};\theta_{1})}{\partial\theta_{1}}\\
		& +\frac{penalty}{n_{inter}}\sum_{i=1}^{n_{inter}}2(u_{1}(\boldsymbol{x}_{i};\theta_{1})-u_{2}(\boldsymbol{x}_{j};\theta_{2}))\frac{\partial u_{1}(\boldsymbol{x}_{i};\theta_{1})}{\partial\theta_{1}}.
	\end{split}
	\label{eq:The derivative of the potential energy function with respect to parameter 1}
\end{equation}
If we consider the one dimension elastic problem, the strain energy
density is
\begin{equation}
	w_{\varepsilon}(u(\boldsymbol{x}_{i};\boldsymbol{\theta}))=\frac{1}{2}E(\frac{\partial u(\boldsymbol{x};\boldsymbol{\theta})}{\partial\boldsymbol{x}}|_{\boldsymbol{x}=\boldsymbol{x}_{i}})^{2}.
\end{equation}
When substituting the strain energy density to \Cref{eq:The derivative of the potential energy function with respect to parameter 1},
we can obtain

\begin{equation}
	\begin{split}\frac{\partial\mathcal{L}}{\partial\theta_{1}} & =\sum_{i=1}^{n_{1}}w_{i}E\frac{\partial u_{1}(\boldsymbol{x};\theta_{1})}{\partial\boldsymbol{x}}|_{\boldsymbol{x}=\boldsymbol{x}_{i}}\frac{\partial(\frac{\partial u_{1}(\boldsymbol{x};\theta_{1})}{\partial\boldsymbol{x}}|_{\boldsymbol{x}=\boldsymbol{x}_{i}})}{\partial\theta_{1}}-\sum_{i=1}^{n_{1}}w_{i}f_{i}\frac{\partial u_{1}(\boldsymbol{x}_{i};\theta_{1})}{\partial\theta_{1}}\\
		& +\frac{penalty}{n_{inter}}\sum_{i=1}^{n_{inter}}2(u_{1}(\boldsymbol{x}_{i};\theta_{1})-u_{2}(\boldsymbol{x}_{j};\theta_{2}))\frac{\partial u_{1}(\boldsymbol{x}_{i};\theta_{1})}{\partial\theta_{1}}.
	\end{split}
	\label{eq:energy to theta}
\end{equation}
Here we use the chain rule to get
\begin{equation}
	\frac{\partial u_{1}(\boldsymbol{x};\theta_{1})}{\partial\boldsymbol{x}}|_{x=x_{i}}=(\frac{\partial u_{1}}{\partial\boldsymbol{z}^{(L+1)}}\frac{\partial\boldsymbol{z}^{(L+1)}}{\partial\boldsymbol{a}^{(L)}})(\frac{\partial\boldsymbol{a}^{(L)}}{\partial\boldsymbol{z}^{(L)}}\frac{\partial\boldsymbol{z}^{(L)}}{\partial\boldsymbol{a}^{(L-1)}})...(\frac{\partial\boldsymbol{a}^{(2)}}{\partial\boldsymbol{z}^{(2)}}\frac{\partial\boldsymbol{z}^{(2)}}{\partial\boldsymbol{a}^{(1)}})(\frac{\partial\boldsymbol{a}^{(1)}}{\partial\boldsymbol{z}^{(1)}}\frac{\partial\boldsymbol{z}^{(1)}}{\partial\boldsymbol{x}})\label{eq:udx},
\end{equation}
where $a^{(M)}$ and $z^{(M)}$ is an activation neuron(act by the
activation function) and a linear neuron (not through an activation
function), $M=1,2,...L$, L is the hidden layer of the neural network.
The derivative term is
\begin{equation}
	\frac{\partial\boldsymbol{a}^{(L)}}{\partial\boldsymbol{z}^{(L)}}=diag(\boldsymbol{\sigma}^{'(L)}|_{\boldsymbol{z}^{(L)}}),
\end{equation}
where $diag(\sigma^{'(L)}|_{z^{(L)}})$ denotes that the diagonal
element is $\sigma^{'(L)}|z_{i}^{(L)}$, $i=1,2,...L_{n}$, $L_{n}$is
the number of the neuron in layer L. Another derivative term is
\begin{equation}
	\frac{\partial\boldsymbol{z}^{(L+1)}}{\partial\boldsymbol{a}^{(L)}}=\boldsymbol{W}^{(L+1)},
\end{equation}
where $W$ is the weight of the layer $L+1$, We obtain a more detailed
expression of \Cref{eq:udx} 

\begin{equation}
	\frac{\partial u_{1}(x;\theta_{1})}{\partial\boldsymbol{x}}|_{\boldsymbol{x}=\boldsymbol{x}_{i}}=[\boldsymbol{W}^{(L+1)}][diag(\boldsymbol{\sigma}^{'(L)}|_{\boldsymbol{z}^{(L)}})\cdot\boldsymbol{W}^{(L)}]...[diag(\boldsymbol{\sigma}^{'(2)}|_{\boldsymbol{z}^{(2)}})\cdot\boldsymbol{W}^{(2)}][diag(\boldsymbol{\sigma}^{'(1)}|_{\boldsymbol{z}^{(1)}})\cdot\boldsymbol{W}^{(1)}].
\end{equation}

\begin{equation}
	\frac{\partial u_{1}(\boldsymbol{x};\theta_{1})}{\partial\theta_{1}^{(r)}}|_{\boldsymbol{x}=\boldsymbol{x}_{i}}=[\boldsymbol{W}^{(L+1)}][diag(\boldsymbol{\sigma}^{'(L)}|_{\boldsymbol{z}^{(L)}})\cdot\boldsymbol{W}^{(L)}]...[diag(\boldsymbol{\sigma}^{'(r)}|_{\boldsymbol{z}^{(r)}})\cdot\boldsymbol{I}_{r_{n}}\bigotimes\boldsymbol{a}^{(r-1)}]
\end{equation}
It is worth noting that, in general, the last layer has no activation
function, so $diag(\sigma^{'(L+1)}|_{z^{(L+1)}})=1$. So the term
in \ref{eq:energy to theta} is 
\begin{equation}
	\begin{split}\frac{\partial(\frac{\partial u(x;\theta)}{\partial\boldsymbol{x}}|_{\boldsymbol{x}=\boldsymbol{x}_{i}})}{\partial\theta_{w}^{(r)}} & =[\boldsymbol{W}^{(L+1)}][diag(\boldsymbol{\sigma}^{'(L)}|_{\boldsymbol{z}^{(L)}})\cdot\boldsymbol{W}^{(L)}]\\
		& ...[diag(\boldsymbol{\sigma}^{'(r)}|_{\boldsymbol{z}^{(L)}})]\bigotimes[diag(\boldsymbol{\sigma}^{'(r-1)}|_{\boldsymbol{z}^{(r-1)}})\cdot\boldsymbol{W}^{(r-1)}]\\
		& ...[diag(\boldsymbol{\sigma}^{'(2)}|_{\boldsymbol{z}^{(2)}})\cdot W^{(2)}][diag(\boldsymbol{\sigma}^{'(1)}|_{\boldsymbol{z}^{(1)}})\cdot\boldsymbol{W}^{(1)}].
	\end{split}
\end{equation}
where $\bigotimes$ denotes the merge action the tensor, so the derivative
of the loss function with respect to the parameters can be expressed
as

\begin{equation}
	\begin{split}\frac{\partial\mathcal{L}}{\partial\theta_{w}^{(r)}} & =\sum_{i=1}^{n_{1}}w_{i}E\{[\boldsymbol{W}^{(L+1)}][diag(\boldsymbol{\sigma}^{'(L)}|_{\boldsymbol{z}^{(L)}})\cdot\boldsymbol{W}^{(L)}]...[diag(\boldsymbol{\sigma}^{'(2)}|_{\boldsymbol{z}^{(2)}})\cdot\boldsymbol{W}^{(2)}][diag(\boldsymbol{\sigma}^{'(1)}|_{\boldsymbol{z}^{(1)}})\cdot\boldsymbol{W}^{(1)}]\}|_{x=x_{i}}\\
		& \cdot\{[\boldsymbol{W}^{(L+1)}][diag(\boldsymbol{\sigma}^{'(L)}|_{\boldsymbol{z}^{(L)}})\cdot\boldsymbol{W}^{(L)}]...[diag(\boldsymbol{\sigma}^{'(r)}|_{\boldsymbol{z}^{(L)}})]\bigotimes[diag(\boldsymbol{\sigma}^{'(r-1)}|_{\boldsymbol{z}^{(r-1)}})\cdot\boldsymbol{W}^{(r-1)}]\\
		& ...[diag(\boldsymbol{\sigma}^{'(2)}|_{\boldsymbol{z}^{(2)}})\cdot\boldsymbol{W}^{(2)}][diag(\boldsymbol{\sigma}^{'(1)}|_{\boldsymbol{z}^{(1)}})\cdot\boldsymbol{W}^{(1)}]\}|_{\boldsymbol{x}=\boldsymbol{x}_{i}}-\sum_{i=1}^{n_{1}}w_{i}f_{i}\{[\boldsymbol{W}^{(L+1)}][diag(\boldsymbol{\sigma}^{'(L)}|_{\boldsymbol{z}^{(L)}})\cdot\boldsymbol{W}^{(L)}]\\
		& ...[diag(\boldsymbol{\sigma}^{'(r)}|_{\boldsymbol{z}^{(r)}})\cdot\boldsymbol{I}_{r_{n}}\bigotimes\boldsymbol{a}^{(r-1)}]\}|_{\boldsymbol{x}=\boldsymbol{x}_{i}}+\frac{penalty}{n_{inter}}\sum_{i=1}^{n_{inter}}2(u_{1}(\boldsymbol{x}_{i};\theta_{1})-u_{2}(\boldsymbol{x}_{j};\theta_{2}))\\
		& [\boldsymbol{W}^{(L+1)}][diag(\boldsymbol{\sigma}^{'(L)}|_{\boldsymbol{z}^{(L)}})\cdot\boldsymbol{W}^{(L)}]...[diag(\boldsymbol{\sigma}^{'(r)}|_{\boldsymbol{z}^{(r)}})\cdot\boldsymbol{I}_{r_{n}}\bigotimes\boldsymbol{a}^{(r-1)}].
	\end{split}
\end{equation}

The neural network 2 is the same as the neural network 1.

\section{The strong form of the hyperelasticity\label{sec:Appendix C hyper strong}}

In the section,  we derive the strong form of hyperelasticity
with respect to Neo-Hookean.  The strong form is 

\textbf{\begin{equation}
\begin{cases}
\nabla_{X}\cdot P+f=0 & x\in\Omega\\
u=\bar{u} & x\in\partial\Omega^{eb}\\
N\cdot P=\bar{t} & x\in\partial\Omega^{nb}
\end{cases}
\end{equation}}
where $P=(\frac{\partial\Psi}{\partial F})^{T}$ and $\Psi=\frac{1}{2}\lambda(lnJ)^{2}-uln(J)+\frac{1}{2}u(trace(C)-3)$,
 so we can use the chain derivative to get
\begin{equation}
	P=(\frac{\partial\Psi}{\partial J}\frac{\partial J}{\partial\boldsymbol{F}}+\frac{\partial\Psi}{\partial trace(\boldsymbol{C})}\frac{\partial trace(\boldsymbol{C})}{\partial\boldsymbol{C}}\frac{\partial\boldsymbol{C}}{\partial\boldsymbol{F}})^{T}\label{eq:Pk1}
\end{equation}

Further,  we can use the tensor derivative to get
\begin{equation}
	\begin{cases}
		\frac{\partial J}{\partial\boldsymbol{F}}=\boldsymbol{J}\cdot\boldsymbol{F}^{-T}\\
		\frac{\partial trace(\boldsymbol{C})}{\partial\boldsymbol{C}}=\boldsymbol{I}\\
		\frac{\partial C_{ij}}{\partial F_{mn}}=\delta_{in}F_{mj}+F_{mi}\delta_{jn}
	\end{cases}\label{eq:derivative of tensor}
\end{equation}

We substitute \Cref{eq:derivative of tensor} to \Cref{eq:Pk1},  we
can obtain
\begin{equation}
	\boldsymbol{P}=\mu\boldsymbol{F}^{T}+[\lambda ln(J)-\mu]\boldsymbol{F}^{-1}
\end{equation}

where 
\begin{equation}
	\boldsymbol{F}=\frac{\partial\boldsymbol{x}}{\partial\boldsymbol{X}}=\boldsymbol{I}+\frac{\partial\boldsymbol{u}(\boldsymbol{X})}{\partial\boldsymbol{X}}
\end{equation}

Finally,  we get the strong form with respect to the $u$
\begin{equation}
	\begin{cases}
		\nabla_{\boldsymbol{X}}\cdot\{\mu(\boldsymbol{I}+\frac{\partial\boldsymbol{u}(\boldsymbol{X})}{\partial\boldsymbol{X}})^{T}+[\lambda ln(J)-\mu](I+\frac{\partial\boldsymbol{u}(\boldsymbol{X})}{\partial\boldsymbol{X}})^{-1}\}+\boldsymbol{f}=0 & \boldsymbol{x}\in\Omega\\
		\boldsymbol{u}=\boldsymbol{\bar{u}} & \boldsymbol{x}\in\partial\Omega^{eb}\\
		\boldsymbol{N}\cdot\{\mu(\boldsymbol{I}+\frac{\partial\boldsymbol{u}(\boldsymbol{X})}{\partial\boldsymbol{X}})^{T}+[\lambda ln(J)-\mu](\boldsymbol{I}+\frac{\partial\boldsymbol{u}(\boldsymbol{X})}{\partial\boldsymbol{X}})^{-1}\}=\boldsymbol{\bar{t}} & \boldsymbol{x}\in\partial\Omega^{nb}
	\end{cases}
\end{equation}

We can find the strong form is more complicated than the energy form. 

\section{Supplementary code}
The code of this work will be available at \url{https://github.com/yizheng-wang/Research-on-Solving-Partial-Differential-Equations-of-Solid-Mechanics-Based-on-PINN}.

\bibliographystyle{elsarticle-num}
\addcontentsline{toc}{section}{\refname}\bibliography{bibtex}

\end{document}